\documentclass[reqno]{amsart}

\usepackage{latexsym,amsthm,amssymb,amsmath,amsfonts,euscript,amstext}
\usepackage{epsfig}

\begin{document}
\title[Local smoothing effect and existence  ]
{Local smoothing effect and existence for a needle crystal growth problem with anisotropic surface tension }

\author
{Xuming Xie}

\address{Xuming Xie \newline
Department of Mathematics\\
Morgan State University\\
Baltimore, MD 21251, USA}
\email{xuming.xie@morgan.edu}

\thanks{Part of this work was done when the author was in residence at MSRI in Berkeley in the Spring of 2011, the author was grateful for its support and hospitality.}
\keywords{Unsteady needle crystal;  initial value problem; well posedness;
 \hfill\break\indent Sobolev space; smothing effect;
 initial-value problem}

\begin{abstract}
 We study an initial value problem for two-dimensional needle crystal growth with anisotropic surface tension. The initial value problem is derived from the so called one-sided model based on complex variables method.
We then obtain the existence and uniqueness of local solution of the needle crystal problem for any initial interface. Furthermore, we obtain that, on average in time, the solution gains $3/2$ derivative of smoothness in spatial variable compared to the initial data. The continuous dependence on the initial data of the solution map is also established.
\end{abstract}

\maketitle
\numberwithin{equation}{section}
\newtheorem{theorem}{Theorem}[section]
\newtheorem{lemma}[theorem]{Lemma}
\newtheorem{remark}[theorem]{Remark}
\newtheorem{definition}[theorem]{Definition}
\newcommand{\abs}[1]{\lvert#1\rvert}
\newcommand{\norm}[1]{\left\Vert#1\right\Vert}

\section{Introduction}
Dendritic growth is one of the earliest and the most profound scientific problems in the area of interfacial pattern formation. This is not only due to its underlying vital technical importance in the material processing industries but also because dendritic growth represents a fascinating class example of nonlinear phenomena in nonequilibrium systems.  From mathematical point of view, dendrite formation  is a free boundary problem like the Stefan problem \cite{Caff, Friedman, Luckhaus,Escher1, Hadzic}. This problem has been extensively studied based on various models. For review papers and books, we refer to Mullins and Sekerka \cite{Mullins}, Langer \cite{Langer}, Kessler et al \cite{Kessler}, Pelce \cite{Pelce},  Ball et al \cite{Ball} and Davis \cite{Davis} and Xu \cite{Xu}.\\

The growth of a single needle crystal  from an undercooled melt is a simplest example of dendritic growth.  It is well known that a steadily moving front is unstable due to Mullins - Sekerka instability. Mullins-Sekerka instability leads to a wide variety of morphologies including compact shapes and dendrites.
The equations studied in the paper were derived in Kunka et al \cite{Kunka1,Kunka2}. They are based on complex variable method which is a very effective technique in handling two dimensional problems in fluid mechanics such as Hele-Show problem \cite{Howison, Xie1, Xie2}. Based on this equations, Kunka et al \cite{Kunka1,Kunka2} studied the linear theory of localized disturbances and a class exact zero-surface-tension solutions if the initial conditions include only poles. They also studied the singular behavior of unsteady dendritical crystal with surface tension. In those situations, a zero of the conformal map that describes the crystal gives birth to a daughter singularity that moves away from the zero and approaches the interface. However, the analysis in \cite{Kunka1,Kunka2} was asymptotic and numerical but not rigorous.\\

 For steady needle crystal, Xie \cite{Xie3,Xie4} proved that 
in the limit of zero surface tension, these
equations  do not have any physically
acceptable solutions when crystalline anisotropy is ignored even though the equations
 admit solutions (Ivantsov solutions) when surface tension is zero. A discrete set of solution was found to  exist when crystalline anisotropy is included. Linear stability of steady needle crystal was also studied \cite{Xie5, Xie7}. Xie \cite{Xie6} established local existence and uniqueness of analytic solution for unsteady crystal with zero surface tension if the initial data is analytic.\\

There have been a variety of literature for similar problems such as Stefan problem and Hele-Shaw problem.
Existence of global weak solutions were established in Chen et al \cite{Chen,chen1}, Luckhaus \cite{Luckhaus}, Kim \cite{Kim,Kim1} and  Almgren and Wang \cite{Almgren}. Duchon \& Robert \cite{Duchon} have proven the existence for short time. Constantin and Pugh \cite{Const} proved the global existence for Hele-Shaw problem with small analytic initial data; Escher and Simonett \cite{Escher1,Escher} discussed solutions in Holder space while Prokert \cite{Prokert1} obtained existence in Sobolev spaces. Gunther and Prockert \cite{Gunther1} obtained existence results for a similar problem with variable surface energy; Friedman and Reitich \cite{Friedman} and Hadzic and Guo \cite{Hadzic} studied the stability of the Stefan problem. Very recently, Xie \cite{Xie8} established local existence and smoothing effect for the Hele-Shaw problem with small initial data.\\

It should also be noted that there were numerous works on related directional solidification 
problems based on ``phase-field models'' and "sharp interface model"; for example,  see J. Langer [62],
 G. Caginalp \cite{Caginalp}, McFadden et al \cite{McFadden} and references therein.\\

 In this paper, we study the unsteady needle crystal problem with anisotropic surface tension in Sobolev space. We obtain the existence and uniqueness of local solution of the needle crystal problem for any initial interface. Furthermore, we obtain that, on average in time, the solution gains $3/2$ derivative of smoothness in spatial variable compared to the initial data. The Lipschitz continuity of the solution map is also established. It is to be noted that unlike some of the afore mentioned results for Stefan or Hele-Shaw problem, our results require  neither the initial data is small nor the surface energy is isotropic.\\

The paper is organized as follows: In this section, we give a description of mathematical formulation of the unsteady needle crystal based on complex variable methods, then we derive an initial value problem  and present the main theorem on the existence and uniqueness of solution. In section 2 we give a few lemmas that will be used in later sections. In section 3, we obtain several existence and apriori estimate results for some linear equations. In particular, we obtain an energy estimates for a regularized linear equation with variable coefficients (see equation (\ref{3.9}) in Lemma 3.3), which is essential to our iteration scheme for the quasilinear equation in section 4. The proof of our main theorem is presented in section 4.

\subsection{Mathematical formulation}

We are interested in the problem of a free dendrite growing in its 
undercooled melt. Temperature is measured in units of $L/c_p$, where L is the
latent heat, $c_p$ the heat capacity. Lengths are measured in units of the tip
 radius of curvature $a$ for the Ivantsov parabola. Velocity is measured in units of U, where U is such that  the Peclet number P defined by
$P=\frac{Ua}{2D}$, where D is the thermal diffusivity. Time is measured in units of $\frac{a}{U}$. The dimensionless undercooling defined as 
$\Delta=(c_p/L)(T_m-T_\infty)$ satisfy Ivantsov relation ($\Delta = \sqrt{\pi P}e^P ~\text{erfc}~ [\sqrt{P}]$),  $T_m$ is the melting temperature and $T_\infty$ is 
the specific temperature an infinity.\\

If heat diffusion in the solid phase is neglected , in the frame where an Ivantsov 
parabolic interface would have been stationary and the crystal-frame coordinates $x-y$ is placed at the tip of the needle crystal, then  the dimensionless temperature T, 
with the melting temperature subtracted,  satisfies 
\begin{equation}
\label{1.1}
2P \frac{\partial T}{\partial t} = 2P\frac{\partial T}{\partial y}+\nabla^2 T
\end{equation}

The condition at infinity that determines T for a specified undercooling is
\begin{equation}
\label{1.2}
T\to-\Delta ~ ~as ~~ y\to \infty
\end{equation}

The conservation of heat through the interface requires
\begin{equation}
\label{1.3}
\frac{\partial T}{\partial n}= -2P (v_n+ \cos \theta )
\end{equation}

where $v_n$ is the normal component of the interface motion and $\theta$ is the angle
 between the interface and y axis. See Fig 1.

\begin{figure}
\begin{center}
\resizebox{0.75\textwidth}{!}
{\includegraphics{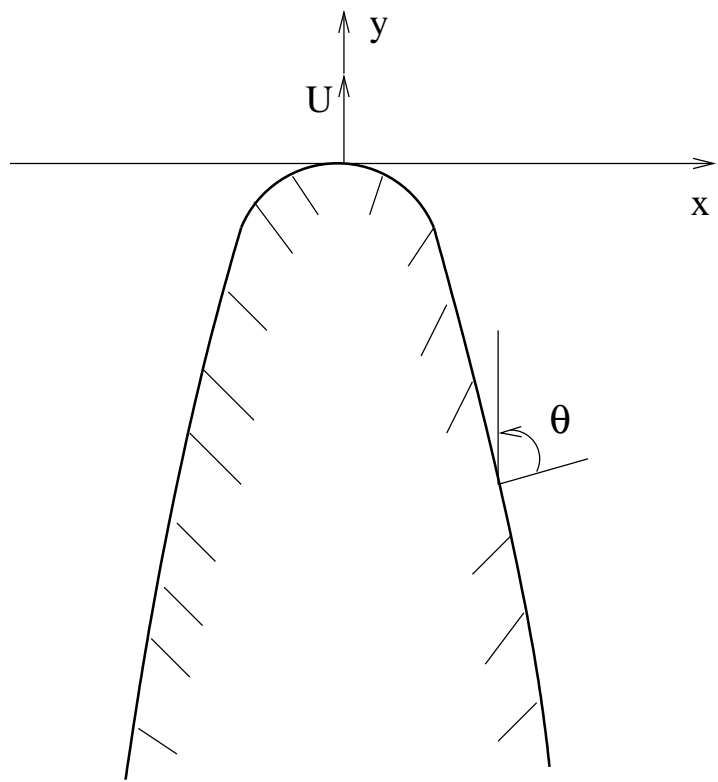}}
\end{center}
\caption{Problem domain: $z= x+iy$ for needle crystal.}
\label{fig:1}
\end{figure}

The Gibbs-Thompson boundary condition at the interface is 
\begin{equation}
\label{1.4}
T=-d_0[1-\gamma\cos 4\theta ]\kappa
\end{equation}
where $d_0$ is a nondimensional capillary parameter defined by
$$
d_0=\frac{dc_p}{aL}T_m
$$
where $d$ is the standard capillary length.\\
In (\ref{1.4}), $\kappa$ is the curvature, $[1-\gamma\cos 4\theta ]$ is included to
 model a standard fourfold anisotropy in the surface energy; $\gamma$ is the crystalline anisotropy.

We consider the conformal  map $z(\xi,t)$ with $\xi=w+is$  that maps the
 upper-half $\xi $ plane into the exterior of the crystal in the z plane, 
where $z=x+iy$. The real $\xi$ axis $s =0$  corresponds the unknown interface (see Fig 2).
 It is clear that determination of the function $z(\xi,t)$ yields the unknown interface. \\

\begin{figure}
\begin{center}
\resizebox{0.75\textwidth}{!}
{\includegraphics{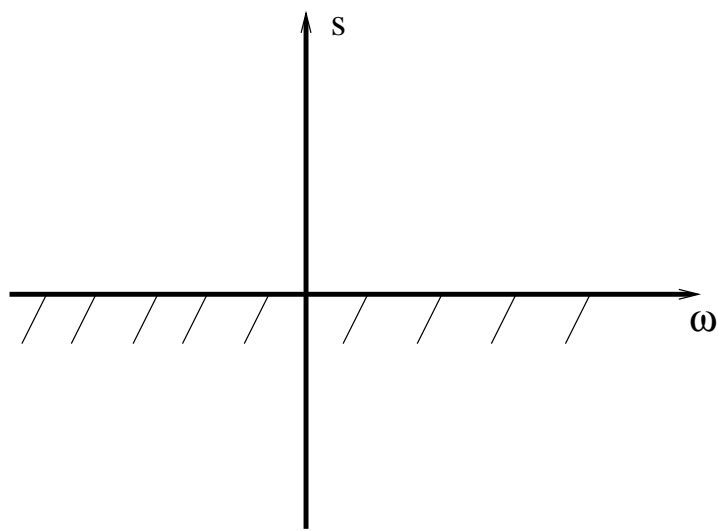}}
\end{center}
\caption{$\xi$ plane: $\xi= w+is$ .}
\label{fig:1}
\end{figure}

Under this transformation, (\ref{1.1}) becomes
\begin{multline}
\label{1.5}
2P|z_\xi|^2\left[ \frac{\partial T}{\partial t}-Re~\left(\frac{z_t}{z_\xi}\right)
\frac{\partial T}{\partial w}-Im~\left(\frac{z_t}{z_\xi}\right)
\frac{\partial T}{\partial s}\right] \\
=
2P\left( Im~(z_\xi)\frac{\partial T}{\partial w}+ 
Re (z_\xi )\frac{\partial T}{\partial s}\right)+\nabla ^2 T .
\end{multline}
(\ref{1.2}) becomes
\begin{equation}
\label{1.6}
T\to -\Delta ~~ as ~~ s\to \infty.
\end{equation}
(\ref{1.3}) becomes 
\begin{equation}
\label{1.7}
\frac{\partial T}{\partial s}= -2P|z_\xi|^2~Im~\left(\frac{z_t+i}{z_\xi}\right).
\end{equation}
The Gibbs- Thompson boundary condition on $s =0$ is 
\begin{equation}
\label{1.8}
T=-d_0[1-\gamma\cos 4\theta ]\kappa(\xi),
\end{equation}
where
\begin{equation}
\label{1.9}
\kappa=-Im~\left[\frac{z_{\xi\xi}}{z_\xi}\right]/|z_\xi |,
\end{equation}
\begin{equation}
\label{1.10}
\cos 4\theta = Re~ \left[\frac{z_\xi^4}{|z_\xi|^4}\right].
\end{equation}
The Ivantsov steady solution corresponds to 
\begin{equation}
\label{1.11}
z_I(\xi)=-i\xi^2 /2 +\xi
\end{equation}
and 
\begin{equation}
\label{1.12}
T_I=-\Delta+\sqrt{\pi P}e^P ~\text{erfc}~[\sqrt{P}(1+s)],
\end{equation}
where
\begin{equation}
\label{1.13}
\Delta = \sqrt{\pi P}e^P ~\text{erfc}~ [\sqrt{P}].
\end{equation}
Now consistent with most experimental conditions, we assume 
\begin{equation}
\label{1.14}
d_0=2\tau P, ~~\epsilon^2=O(1) ~~ as ~~ P\to 0.
\end{equation}
We then use regular perturbation expansion
\begin{equation}
\label{1.15}
T=P T_0 + O(P^2).
\end{equation}
Then , to O(P), (\ref{1.5}) becomes:
\begin{equation}
\label{1.16}
\nabla^2 T_0=0,
\end{equation}
\begin{equation}
\label{1.17}
\frac{\partial T_0}{\partial s}= -2|z_\xi|^2~Im~\left(\frac{z_t+i}{z_\xi}\right),
\end{equation}
\begin{equation}
\label{1.18}
T_0=-2\tau [1-\gamma \cos 4\theta ]\kappa(\xi).
\end{equation}
Since $T_0$ is a harmonic function in two dimension, we can define an analytic 
function $\tilde{W}$ so that 
\begin{equation}
\label{1.19}
T_0=Re (\tilde{W}).
\end{equation}
we decompose $\tilde{W}$
\begin{equation}
\label{1.20}
\tilde{W}=2i\xi -2 W(\xi, t).
\end{equation}
Then (\ref{1.17}) implies that on $s=0$,
\begin{equation}
\label{1.21}
1-Im (W_\xi(\xi ,t))=|z_\xi|^2 Im~\left[ \frac{(z_t+i)}{(z_\xi}
\right],
\end{equation}
and (\ref{1.18}) becomes
\begin{equation}
\label{1.22}
Re~ W(\xi, t)=\tau[1-\gamma \cos 4\theta]\kappa(\xi)\text{ on } s=0.
\end{equation}

For much smaller times, it is appropriate to invoke the Ivantsov solution behavior (\ref{1.12}) and assume that in the far field
\begin{equation}\label{1.23}
\begin{split}
&T_0\to -2s   \text{ as } s\to \infty,\\
&z(\xi,t)-z^I(\xi)\to 0,~~z_\xi(\xi,t)-(z_I)_\xi(\xi)\to 0 \text{ as $\xi$ real and } \xi \to \infty;
\end{split}
\end{equation}
which implies
\begin{equation}\label{1.24}
Re~W(\xi,t)\to 0 \text{ as }\xi\to \infty.
\end{equation}
Since $W(\xi,t)$ is analytic in the upper half plane, from Plemelj's formula (see \cite{Carrier}) and (\ref{1.22}), we obtain
\begin{equation}\label{1.25}
W(\xi,t)=\frac{1}{\pi i} (P)\int_{-\infty}^\infty \frac{\tau[1-\gamma \cos 4\theta]\kappa(\xi')d\xi'}{\xi'-\xi}
\end{equation}
(\ref{1.23}) implies that $Im~\left(\frac{z_t+i}{z_\xi}\right) \to 0$ as $\xi\to\infty$.
Since $z(\xi,t)$ is analytic in the upper half plane and $z_\xi\ne 0$, from Plemelj's formula and (\ref{1.21}), we obtain in the upper half plane $s=Im~\xi>0$
\begin{equation}\label{1.26}
\frac{z_t+i}{z_\xi}=\frac{1}{\pi } \int_{-\infty}^\infty \frac{[1-Im(W_\xi)]d\xi'}{|z_\xi|^2(\xi'-\xi)}.
\end{equation}
Letting $s\to 0^+$ in (\ref{1.26}), and using (\ref{1.25}), we obtain
\begin{equation}\label{1.27}
z_t+i=z_\xi [H[Q](\xi,t)+iQ[z](\xi,t)],~ \text{for real }\xi,
\end{equation}
where $H$ is the Hilbert transform
\begin{equation}\label{1.28}
H[Q](\xi,t)=\frac{1}{\pi}(P)\int_{-\infty}^{\infty}\frac{Q[z](\xi',t)}{(\xi'-\xi)}d\xi',
\end{equation}
and
\begin{equation}\label{1.29}
Q[z](\xi,t)=\frac{1+\tau\partial_\xi H[(1-\gamma \cos 4\theta )\kappa ]}{|z_\xi|^2}.
\end{equation}
The initial condition is
\begin{equation}\label{1.30}
z(\xi,0)=z_0(\xi).
\end{equation}

We introduce
\begin{equation}\label{1.31}
\log z_\xi=\ln |z_\xi|+i \arg z_\xi=h(\xi)+iq(\xi).
\end{equation}
From (\ref{1.31}), we have $h_\xi+iq_\xi=\frac{z_{\xi\xi}}{z_\xi}$ is analytic in the upper half plane. (\ref{1.22}) implies that $h_\xi,q_\xi\to 0$ as $\xi\to \infty$, from Plemelj formula
\begin{equation}\label{1.32}
q_\xi(\xi,t)=-H[h_\xi(\xi,t)], \text{ for real} \xi.
\end{equation}
then  we have from (\ref{1.9}),(\ref{1.31}) and (\ref{1.32})
\begin{equation}\label{1.33}
\kappa=\frac{H[h_\xi]}{e^h},
\end{equation}
and  from (\ref{1.10})
\begin{equation}\label{1.34}
\cos 4\theta=\cos 4q ;
\end{equation}
so
\begin{equation}\label{1.35}
Q[h]=\frac{1+\tau \partial_\xi H[ (1-\gamma \cos 4q)e^{-h}H[h_\xi]]}{e^{2h}}.
\end{equation}
Taking derivative with respect to $\xi$  in equation (\ref{1.27}) and using (\ref{1.31}), we obtain that $h$ satisfies
\begin{equation}\label{1.36}
h_t=h_\xi H[Q[h]]-q_\xi Q[h]+\partial_\xi\left(H[Q[h]]]\right)
\end{equation}

Using (\ref{1.11}), the corresponding Ivantsov solution can be written  as $\log (z_I)_\xi=h_I+iq_I$, where
\begin{equation}\label{1.37}
h^I(\xi)=\frac{1}{2}\ln (1+\xi^2), q^I(\xi)=-\arcsin \frac{\xi}{\sqrt{1+\xi^2}}.
\end{equation}
We decompose 
\begin{equation}\label{1.38}
u(\xi,t)=h(\xi,t)-h^I(\xi),~v(\xi,t)=q(\xi,t)-q^I(\xi).
\end{equation}
From (\ref{1.31}) and(\ref{1.37}), we have that $u+iv=\log \frac{z_\xi}{(z_I)_\xi}$ is analytic in the upper half plane. (\ref{1.23}) implies that$u,v\to 0$ as $\xi\to \infty$, from Plemelj formula
\begin{equation}\label{1.39}
v(\xi,t)=-H[u(\xi,t)], \text{for real } \xi.
\end{equation}
From (\ref{1.36}) $ u(\xi,t)$ satisfies
\begin{equation}\label{1.40}
u_t=(u+h^I_\xi) H[Q[u]]-(q^I_\xi-H[u]) Q[u] +\partial_\xi\left(H[Q[u]\right)
\end{equation}
where
\begin{equation}\label{1.41}
Q[u]=\frac{1+\tau \partial_\xi H\left[\left(1-\gamma \cos 4(q^I-H[u])\right)e^{-(h^I+u)}H[h^I_\xi+u_\xi]\right]}{e^{2(h^I+u)}}.
\end{equation}
The initial condition is
\begin{equation}\label{1.42}
u(\xi,0)=u_0(\xi)
\end{equation}

\subsection{Main results}

In this paper, $H^s(R)$ denotes the Sobolev space over the real line $R$, with the norm $\norm{f}_{H^s(R)}=(\int (1+|\lambda|)^{2s}|\hat{f}(\lambda)|^2d\lambda )^{1/2}$, where $\hat{f}(\lambda)=\frac{1}{\sqrt{2\pi}}\int f(\xi)e^{-i\xi\lambda}d\xi$ is the Fourier transform of $f(\xi)$. If $s$ is a nonnegative integer $k$, we also use $\norm{f}_{H^k(R)}=\sum_{j=0}^k\norm{\partial_\xi^j f}_{L^2(R)}$. $P(D)$  denotes the pseudodifferential operator defined by
\begin{equation*}
P(D)f(\xi,t)=\frac{1}{\sqrt{2\pi}} \int P(\lambda )\hat{f}(\lambda)e^{i\xi\lambda}d\lambda.
\end{equation*}
Note that $\hat{H(u)}(\lambda )=i~\text{sgn }(\lambda ) \hat{u}$, where $\text{ sgn }(\lambda) =1 $ if $\lambda >0$ and $\text{ sgn }(\lambda) =-1 $ if $\lambda <0$; and $|D|=-H\partial_\xi =-\partial_\xi H$.\\

We are going to prove the main theorem in this paper:
\begin{theorem}
Assume that $\beta=\tau (1-\gamma)>0, u_0\in H^{s+\frac{1}{2}}(R), s\ge 5$ is an integer, $M_0=\norm{u_0}_{H^{s+\frac{1}{2}}(R)}$, then there is  $T=T(M_0)>0$ such that the initial value problem (\ref{1.40}) (\ref{1.42}) has a unique solution $u\in C([0,T], H^{s+\frac{1}{2}}(R))\cap C^1([0,T],H^{s-5/2}(R))\cap L^2([0,T],H^{s+2}(R))$
satisfying
\begin{equation}\label{1.43}
\int_0^T\int_0^{2\pi} |\partial_\xi^{s+2} u(\xi,t)|^2d\xi dt\le CM_0^2.
\end{equation}
Moreover, the solution map $u_0\to u(t)$ from $H^{s+\frac{1}{2}}(R) $ to $ C([0,T], H^{s+\frac{1}{2}}(R) )$\\$\cap L^2([0,T],H^{s+2} (R))$ is Lipschitz continuous.
\end{theorem}
The main idea of the proof is outlined now.  We write the equation (\ref{1.40}) as a quasilinear equation
\begin{equation}\label{1.44}
u_t+B[u] H[\partial^3_\xi u]=\mathcal{N}(u) ;
\end{equation}
where $B[u]\ge \beta$, $\mathcal{N}(u)$ consists of lower order terms compared to $ H[\partial^3_\xi u]$. More precisely, we will show (see Lemma 4.3 and Lemma 4.6))
\begin{equation}\label{1.45}
\norm{\mathcal{N}(u)}_{H^{s-2}}\leq C(\norm{u}_{H^s}).
\end{equation}
 The corresponding linear equation is 
\begin{equation}\label{1.46}
u_t+ b(\xi,t) H[\partial^3_\xi u]=f(\xi);
\end{equation}
where $b(\xi,t)\ge \beta>0$.
It can be shown  ( Lemma 3.3 ) that the solution of the initial value problem (\ref{1.46}) and (\ref{1.42}) has local smoothing effect (Lemma 3.3).  To establish existence theory for linear equations, we use  viscosity method by introducing a regularized term $\epsilon \partial_\xi^6 u$. Finally, an iteration scheme then enable us to carry on the smoothing effect to the quasilinear equation (\ref{1.44}).\\

\section{Preliminary lemmas}

We  introduce some preliminary lemmas.
\begin{lemma}
If $u,v\in H^s(R),$  for any $s>\frac{1}{2}$, then $uv \in H^s(R)$ and
\begin{equation}\label{2.1}
\norm{uv}_{H^s(R)}\leq C \norm{u}_{H^s(R)}\norm{v}_{H^s(R)}.
\end{equation}
\end{lemma}
\begin{lemma}
If $f: R\to R$ is a smooth function with $f(0)=0$, and $u\in H^s(R),\norm{u}_{H^s(R)}\leq M$  for any $s>\frac{1}{2}$, then $f(u)\in H^s(R)$  and
\begin{equation}\label{2.2}
\norm{f(u)}_{H^s(R)}\leq C\norm{u}_{H^s(R)}
\end{equation}
where $C$ depends only on $f, s$ and $M$.
\end{lemma}
The proof of these lemmas can be found in \cite{Taylor}.
\begin{lemma}
(1) For any $u\in H^{\frac{1}{4}}(R), \norm{u}_{L^{4}(R)}\leq K_0\norm{|D|^{1/4}u}_{L^2(R)}$.\\
(2) Let $f\in H^1(R)$ and $g\in H^{1/2}(R)$, then $fg\in H^{1/2}(R)$ and $\norm{fg}_{H^{1/2}(R)}\leq K_0\norm{f}_{H^1(R)}\norm{g}_{H^{1/2}(R)}$.
\end{lemma}
\begin{proof}
See Lemma 7.1 in \cite{Wu}.
\end{proof}
\begin{lemma}
Let $u\in H^s(R)$, then $H(u)\in H^s(R)$ and $\norm{H(u)}_{H^s(R)}\leq \norm{u}_{H^s(R)}$.
\end{lemma}
\begin{proof}
The Lemma follows from the fact  $\hat{H(u)}(\lambda )=i~\text{sgn }(\lambda ) \hat{u}$, where $\text{ sgn }(\lambda) =1 $ if $\lambda >0$ and $\text{ sgn }(\lambda) =-1 $ if $\lambda <0$.
\end{proof}

The Sobolev embedding theorem implies 
\begin{lemma}
Let $u\in H^{r} $, $r>\frac{1}{2}$, then $u\in C(-\infty,\infty) $   and $\norm{u}_{L^\infty  }\leq K_1 \norm{u}_{H^{r} }$, $\norm{u}_{H^p }\leq \norm{u}_{H^{r} }$ for $0\le p<r$.
\end{lemma}
\begin{lemma}
Let $0\le r<p$. Then for every $\eta >0$, there exists a $K(\eta)$ such that
$\norm{u}_{H^r}\le \eta \norm{u}_{H^p}+K(\eta)\norm{u}_{L^2}$.
\end{lemma}
\begin{proof}
\begin{equation*}
\begin{split}
\norm{u}^2_{H^r}&\leq \int_{-\infty}^{\infty}|\lambda|^{2r}|\hat{u}(\lambda)|^2d\lambda+\int_{-\infty}^\infty \hat{u}^2(\lambda)d\lambda\\
&\le \int_{|\lambda|\le N}|\lambda|^{2r}|\hat{u}(\lambda)|^2+N^{2(r-p)}\int_{|\lambda|\ge N}|\lambda|^{2p}|\hat{u}(\lambda)|^2+\int_{-\infty}^\infty \hat{u}^2(\lambda)d\lambda\\
&\le (N^{2r}+1)\norm{u}_{L^2}^2+N^{2(r-p)}\norm{u}_{H^p}^2
\end{split}
\end{equation*}
Since $r-p<0$, the Lemma follows from this if $N$ is chosen large enough.
\end{proof}
\begin{lemma}
Let $f, g \in L^2(R)$, then\\

(1) $\int H[f]H[g]d\xi=\int fg d\xi, ~~~\int gH[f]d\xi =-\int fH[g]d\xi$.\\

(2) $H[H[f]]=-f$.\\

(3) 
\begin{equation} \label{2.3}
H[fg]=H[H[f]H[g]]+gH[f]+fH[g].
\end{equation}
\end{lemma}
\begin{proof}
(1) 
\begin{equation*}
\int H[f]H[g]d\xi=\int \hat{H[f]}(\hat{H[g]})^*d\lambda=\int \hat{f}(\hat{g})^* d\lambda=\int fg d\xi. 
\end{equation*}
where $*$ denotes complex conjugate.\\
\begin{multline*}
\int gH[f]d\xi=\int \hat{H[f]}(\hat{g})^*d\lambda=\int i\text{ sgn }(\lambda)\hat{f}(\hat{g})^* d\lambda\\
= -\int \hat{f}(\hat{H[g]})^*d\lambda=-\int fH[g]d\xi
\end{multline*}
(2) From Plemelj's formula, it is well known that $u(\xi)+iv(\xi)$ can be analytically extended to the upper half complex plane if and only if $v(\xi)=-H[u]$. Since $i(u-iH[u])=H[u]+iu$ can be analytically extended to the upper half complex plane, we have $u=-H[H[u]]$.\\
(3) Since $(f-iH[f])(g-iH[g])=fg-H[f]H[g]-i(fH[g]+gH[f])$ can be analytically extended to the upper half complex plane, we have
\begin{equation*}
(fH[g]+gH[f])=H[fg-H[f]H[g]]
\end{equation*}
which implies (\ref{2.3})
\end{proof}
We define commutator $[H,f]g=H[fg]-fH[g]$, and we have
\begin{lemma}
If $f\in H^s$, $s\ge 3$, then $\norm{[H,f]g}_{H^s}\leq C \norm{g}_{H^{s-2}}\norm{f}_{H^s}$.
\end{lemma}
\begin{proof}
See Corollary 3.8 in \cite{Ambrose}.
\end{proof}

\section{Linear equations}
In this section, we obtain some results for some linear equations. In this and following sections, $C>0$ represents a generic constant, it may vary from line to line.
\begin{lemma} Let $f\in L^2([0,T],H^{s-3}(R)), s\ge 3$, and $u_0\in H^s(R)$, then the following initial value problem:
\begin{equation}\label{3.1}
u_t-\epsilon \partial^6_\xi u=f(t,\xi), u|_{t=0}=u_0(\xi).
\end{equation}
has an unique solution $u(t,\xi)\in C([0,T],H^s(R))\cap L^2([0,T],H^{s+3}(R))$ and for every $t$ and any $\epsilon>0$
\begin{equation}\label{3.2}
\norm{u}^2_{H^s}+c_0\epsilon e^t\int_0^t\norm{u}^2_{H^{s+3}(R)}e^{-\rho}d\rho\leq \frac{e^t}{2c_0\epsilon}\int_0^te^{-\rho}\norm{f}^2_{H^{s-3}(R)}d\rho +e^t\norm{u_0}^2_{H^s(R)}.
\end{equation}
where $c_0>0$ is some constant.
\end{lemma}
\begin{proof}
We first establish the energy inequality for  $u\in C([0,T],H^s(R))$\\
$\cap L^2([0,T],H^{s+3}(R))$ satisfying (\ref{3.1}). Taking Fourier transform in (\ref{3.1}) leads
\begin{equation}\label{3.3}
\hat{u}_t+\epsilon\lambda^6 \hat{u}=\hat{f}(\lambda,t).
\end{equation}
so we have
\begin{equation}\label{3.4}
\begin{split}
\frac{d}{dt}\norm{u}^2_{H^s}&=\partial_t\int (1+|\lambda|)^{2s}\hat{u}^2d\lambda =2\int (1+|\lambda|)^{2s}\hat{u}\hat{u}_td\lambda\\
&=2\int (1+|\lambda|)^{2s}\hat{u}\hat{f}d\lambda -2\epsilon\int (1+|\lambda|)^{2s}\lambda^6\hat{u}^2d\lambda
\end{split}
\end{equation}
Since there exists a positive constant $c_0$ depending on $\epsilon$ such that for all $\lambda \in (-\infty,\infty)$
\begin{equation}\label{3.5}
(1+|\lambda|)^{2s}(\lambda^6 +\frac{1}{2\epsilon})\ge c_0(1+|\lambda|)^{2s+6},
\end{equation}
which  implies
\begin{equation}\label{3.6}
\int (1+|\lambda|)^{2s}\lambda^6\hat{u}^2d\lambda\ge c_0\norm{u}^2_{H^{s+3}}-\frac{1}{2\epsilon}\norm{u}_{H^s}.
\end{equation}
Young's inequality implies
\begin{multline}\label{3.7}
\int (1+|\lambda|)^{2s}\hat{u}\hat{f}d\lambda\le \epsilon c_0\int (1+|\lambda|)^{2s+6}\hat{u}^2d\lambda +\frac{1}{\epsilon 2 c_0}\int (1+|\lambda|)^{2(s-3)}\hat{f}^2d\lambda\\
\leq \epsilon c_0\norm{u}^2_{H^{s+3}(R)}+\frac{1}{2c_0\epsilon}\norm{f}_{H^{s-3}}^2.
\end{multline}
(\ref{3.4}) - (\ref{3.7}) lead to 
\begin{equation}\label{3.8}
\frac{d}{dt}\norm{u}_{H^s}^2\le \norm{u}^2_{H^s}-c_0\epsilon \norm{u}^2_{H^{s+3}}+\frac{1}{2\epsilon c_0}\norm{f}^2_{H^{s-3}},
\end{equation}
which leads to (\ref{3.2}).\\
The existence and uniqueness of solution to ( \ref{3.1}) then follows from the standard linear theory by combining (\ref{3.2}) and a Galerkin approximation.

\end{proof}

We are going to prove the main theorem using an iteration scheme. To the end, we need to study a regularized linear equation with variable coefficients, i.e  the equation (\ref{3.9}) in the following lemma.
\begin{lemma}
Let $\epsilon>0, T>0, s\ge 3\frac{1}{2} $ be constants, $f\in L^2([0,T], H^{s-3}(R))$, $b(\xi,t)\in C([0,T],H^{s-3})$
and $u_0\in H^{s}(R)$. Then the initial value problem
\begin{equation}\label{3.9}
\begin{split}
&u_t+b(\xi,t) H[\partial^3_\xi u]
-\epsilon\partial^6_\xi u=f(t,\xi),\\
 &u|_{t=0}=u_0
\end{split}
\end{equation}
has a unique solution $u\in C([0,T],H^s(R))\cap L^2([0,T],H^{s+3}(R))$.
\end{lemma}
\begin{proof}
Let $u^0(\xi,t)=0$; we construct a sequence of $u^k(\xi,t)\in C([0,T],H^s(R))\cap L^2([0,T],H^{s+3}(R))$ through
\begin{equation}\label{3.10}
\begin{split}
&u^{k+1}_t-\epsilon\partial^6_\xi u^{k+1}
=f(t,\xi)-b H[\partial^3_\xi u^k],\\
& u^{k+1}|_{t=0}=u_0
\end{split}
\end{equation}
Since $f\in L^2([0,T], H^{s-3}(R))$, $b(\xi,t)\in C([0,T],H^{s-3})$ and $u^k\in C([0,T],H^s(R))$, the right hand side of (\ref{3.9}) is in $L^2([0,T],H^{s-3}(R))$, applying Lemma 3.1, we have  $u^{k+1}\in C([0,T],H^s(R))\cap L^2([0,T],H^{s+3}(R))$, and 
\begin{equation}\label{3.11}
\begin{split}
&(u^{k+1}-u^k)_t-\epsilon\partial^6_\xi (u^{k+1}-u^k)
 =-b H[\partial^3_\xi (u^k-u^{k-1})]\\
 &(u^{k+1}-u^k)|_{t=0}=0
\end{split}
\end{equation}
Applying (\ref{3.2}), we have
\begin{equation}\label{3.12}
\begin{split}
&\norm{u^{k+1}-u^k}^2_{H^s}+c_0\epsilon e^t\int_0^t\norm{u^{k+1}-u^k}^2_{H^8(R)}e^{-\rho}d\rho\\
&~ \leq \frac{e^t}{2c_0\epsilon}\int_0^te^{-\rho}\norm{b H(\partial^3_\xi (u^k-u^{k-1}))}_{H^{s-3}(R)}d\rho\\
&~ \leq  \frac{e^t}{2c_0\epsilon}\int_0^te^{-\rho}\norm{b}_{H^{s-3}}\norm{u^k-u^{k-1}}^2_{H^s(R)}d\rho\\
&~ \leq c_1\int_0^t\norm{u^k-u^{k-1}}^2_{H^s(R)}d\rho
\end{split}
\end{equation}
where constant $c_1$ depends on $\norm{b}_{L^\infty([0,T],H^{s-3}(R))}$ and $c_0$. This shows that $\{ u^k\}$ is a Cauchy sequence in $C([0,T],H^s(R))\cap L^2([0,T],H^{s+3}(R))$. Therfore, 
\begin{equation}\label{3.13}
u^k\to u \text{ in } C([0,T],H^s(R))\cap L^2([0,T],H^{s+3}(R)).
\end{equation}
and $u \in C([0,T],H^s(R))\cap L^2([0,T],H^{s+3}(R))$ is the unique solution of (\ref{3.9}).
\end{proof}

The following lemma gives key energy estimates and is essential to prove the main theorem.

\begin{lemma}
Let $T>0,\beta>0, \epsilon\ge 0 $ be constants, $s\ge 3$ be an integer. Let $b\in C([0,T],H^{s}), b\ge\beta$, $f\in L^2([0,T], H^{s-1} )$, 
and if $u\in L^\infty([0,T], H^{s+\frac{1}{2}} )\cap C^1([0,T], H^{r} )\cap L^2([0,T], H^{s+2} )$ ($r$ is any real number) and $\epsilon |D|^{s+3\frac{1}{2}}u\in L^2([0,T]\times R)$ is a solution of the initial value problem (\ref{3.9}),
then $u\in C([0,T], H^{s+\frac{1}{2}} )$ and for every $t\in [0,T]$ and for $j=0,1.$
\begin{multline}\label{3.14}
\norm{u(\cdot,t)}^2_{H^{s-j+1/2} }+\beta \int_0^t e^{K(t-\rho)} \norm{\partial_\xi^{s-j+2}u(\cdot,\rho)}^2_{L^2 }d\rho\\
+\epsilon \int_0^t e^{K(t-\rho)}\norm{u}^2_{H^{s-j+3\frac{1}{2}} }d\rho \\
\leq e^{Kt}\norm{u_0(\cdot,t)}^2_{H^{s-j+1/2} }+K_1\int_0^t e^{K(t-\rho)}\norm{f(\cdot,\rho)}^2_{H^{s-j-1} }d\rho,
\end{multline}
where $K,K_1$ are positive constant independent of $\epsilon$.\\
\end{lemma}

\begin{proof}
Assume that $u\in C^1([0,T],H^{s+9} ),f\in L^2([0,T], H^{s} )$. From equation ( \ref{3.9}), we get
\begin{equation}\label{3.15}
\begin{split}
&\frac{d}{dt}\int\left\vert|D|^{1/2}(\partial^s_\xi u)\right\vert^2 d\xi=2\int (|D|\partial^s_\xi u)(\partial^s_\xi u)_td\xi \\
&=2\int (|D|\partial^s_\xi u)\left[ \partial^s_\xi f+\partial^{s}_\xi\left(bH(\partial_\xi^3 u)\right)+\epsilon\partial^{s+6}_\xi u\right]d\xi
\end{split}
\end{equation}
Now we estimate each term in the right hand side of (\ref{3.15}).\\
Using integration by parts and  Young's inequality, we have
\begin{multline}\label{3.16}
2\int (|D|\partial^s_\xi u) \partial^s_\xi fd\xi =-2\int (|D|\partial^{s+1}_\xi u) \partial^{s-1}_\xi fd\xi\\
\leq \epsilon_0\int (|D|\partial^{s+1}_\xi u)^2d\xi +\frac{C}{\epsilon_0}\int (\partial^{s-1}_\xi f)^2d\xi,
\end{multline}
Using integration by parts and Lemma 2.7, we obtain
\begin{equation}\label{3.17}
\epsilon \int (|D|\partial^s_\xi u)(\partial^{s+6} u)=-\epsilon \int (|D|^{1/2}\partial_\xi ^{s+3}u)^2 d\xi;
\end{equation}
We used Cauchy inequality, integration by parts and Lemma 2.7 .

We write
\begin{multline}\label{3.18}
2\int (|D|\partial^s_\xi u)\partial^{s}_\xi\left(bH[\partial_\xi^3 u]\right)\\
=2\int (|D|\partial^s_\xi u)b|D|\partial^{s+2}_\xi  u
+2\int (|D|\partial^s_\xi u)\sum_{k=1}^{s}(\partial^{k}_\xi b)|D|\partial_\xi^{s-k+2} u.
\end{multline}
By integration by parts, the first term in (\ref{3.18}) can be estimated as
\begin{multline}\label{3.19}
2\int (|D|\partial^s_\xi u)b|D|\partial^{s+2}_\xi  u=-2\int \partial_\xi[(|D|\partial^s_\xi u)b]|D|\partial^{s+1}_\xi  u\\
=-2\int b[|D|\partial^{s+1}_\xi  u]^2-2\int \partial_\xi b(|D|\partial^s_\xi u)|D|\partial^{s+1}_\xi  u\\
\leq -2\beta\int [|D|\partial^{s+1}_\xi  u]^2-2\int \partial_\xi b(|D|\partial^s_\xi u)|D|\partial^{s+1}_\xi  u
\end{multline}
so (\ref{3.18}) becomes
\begin{multline}\label{3.20}
2\int (|D|\partial^s_\xi u)\partial^{s}_\xi[bH(\partial_\xi^3 u)]\\
\leq -2\beta\int [|D|\partial^{s+1}_\xi  u]^2+2\int (|D|\partial^s_\xi u)\sum_{k=2}^{s}(\partial^{k}_\xi b)|D|\partial_\xi^{s-k+2} u.
\end{multline}
The $k=s$ term  in the sum on the right hand side of (\ref{3.20}) can be estimated as follows:
\begin{multline}\label{3.21}
2\int (|D|\partial^s_\xi u)(\partial^s_\xi b)(|D|\partial^{2}_\xi u)d\xi\\
\leq 2\left(\int(\partial^s_\xi b)^2d\xi \right)^{1/2}\left(\int (|D|\partial^s_\xi u)^4d\xi\right)^{1/4} \left(\int (|D|\partial^2_\xi u)^4d\xi\right)^{1/4}\\
\leq 2K_0^2\norm{b}_{H^s}\norm{u}_{H^{s+5/4}}\norm{u}_{H^{3+1/4}}\leq 2K_0^2\norm{b}_{H^s}\norm{u}^2_{H^{s+5/4}}\\
\leq \epsilon_0 \norm{u}^2_{H^{s+2}}+K(\epsilon_0)\norm{u}_{L^2}.
\end{multline}
where $\epsilon_0$ is a small positive constant which will be determined later. Here we used Cauchy's inequality and Lemma 2.3 and Lemma 2.6.\\
The $k=2$ term  in the sum on the right hand side of (\ref{3.20}) can be estimated as follows
\begin{multline}\label{3.22}
\int (|D|\partial^s_\xi u)^2(\partial^2_\xi b)d\xi\leq \left(\int (|D|\partial^s_\xi u)^4d\xi\right)^{1/2} \left(\int(\partial^2_\xi b)^2d\xi \right)^{1/2}\\
\leq 2K_0\norm{b}_{H^s}\norm{u}^2_{H^{s+5/4}}
\leq \epsilon_0 \norm{u}^2_{H^{s+2}}+K(\epsilon_0)\norm{u}_{L^2};
\end{multline}
 Here we used Cauchy's inequality and Lemma 2.3 and Lemma 2.6.\\
The remaining terms in the sum on the right hand side of (\ref{3.20}) can be estimated in the same fashion. So (\ref{3.15}) leads to
\begin{multline}\label{3.23}
\frac{d}{dt}\int \left\vert |D|^{1/2}(\partial_\xi^s u)\right\vert^2 d\xi
 \leq \frac{C}{\epsilon_0}\int (\partial_\xi^{s-1} f)^2 d\xi\\
-\epsilon \int (|D|^{1/2}\partial_\xi ^{s+3}u)^2 d\xi
+((2s-1)\epsilon_0-2\beta)\norm{u}_{H^{s+2}} +K\norm{u}_{H^{s+\frac{1}{2}}}
\end{multline}
where $K$ is a positive constant depending on $\epsilon_0$ and $\norm{b}_{H^s}$.
Now
\begin{equation}\label{3.24}
\begin{split}
&\frac{d}{dt}\int\left\vert|D|^{1/2}(\partial^{s-1}_\xi u)\right\vert^2 d\xi=2\int (|D|\partial^{s-1}_\xi u)(\partial^{s-1}_\xi u)_td\xi \\
&=2\int (|D|\partial^{s-1}_\xi u)\left[ \partial^{s-1}_\xi f+\partial^{s-1}_\xi\left(bH[\partial_\xi^3 u]\right)+\epsilon\partial^{s-1+6}_\xi u\right]d\xi
\end{split}
\end{equation}
In the same fashion, we can estimate each term of (\ref{3.24}) and  we obtain
\begin{multline}\label{3.25}
\frac{d}{dt}\int \left\vert |D|^{1/2}(\partial_\xi^{s-1} u)\right\vert^2 d\xi
 \leq \frac{C}{\epsilon_0}\int (\partial_\xi^{s-2} f)^2 d\xi\\
-\epsilon \int (|D|^{1/2}\partial_\xi ^{s+2}u)^2 d\xi
+((2s-1)\epsilon_0-2\beta)\norm{u}_{H^{s+1}} +K\norm{u}_{H^{s-1+\frac{1}{2}}}
\end{multline}
and
\begin{multline}\label{3.26}
\frac{d}{dt}\int |u|^2 d\xi =2\int uu_td\xi
=2\int u[f+bH[\partial_\xi^3u]+\epsilon \partial^6u_\xi]\\
\leq 2\int uf -\epsilon \int (\partial^3_\xi u)^2 
 +2\int [b_\xi u+b u_\xi]H[\partial_\xi^2u]
\\
\leq \int u^2+\int f^2 -\epsilon \int (\partial^3_\xi u)^2
+[\int(\partial^2_\xi u)^2 +\norm{b}_\infty^2\int u_\xi^2+\norm{b_\xi}_\infty^2 \int u^2]\\
\leq C\int u^2+\int f^2 -\epsilon \int (\partial^3_\xi u)^2+C[\int u_\xi^2+\int (\partial_\xi^2 u)^2],
\end{multline}
Now we can choose $\epsilon_0,\delta_1$ small enough so that in (\ref{3.23}) and (\ref{3.25}) 
$$ ((2s-1)\epsilon_0-2\beta)\leq -\beta$$

Adding (\ref{3.23}), (\ref{3.25}) to (\ref{3.26}), we get for $j=0,1$
\begin{multline}\label{3.27}
\frac{d}{dt}\norm{u}^2_{H^{s-j+\frac{1}{2}}}+\epsilon\norm{u}^2_{H^{s-j+3\frac{1}{2}}}
+\beta \norm{\partial_\xi^{s-j+3}u}^2_{L^{2}}\leq K_1\norm{f}_{H^{s-j-1}}+K\norm{u}_{H^{s-j+\frac{1}{2}}},
\end{multline}
which lead to (\ref{3.14}). So we have proven (\ref{3.14}) for $u\in C^1([0,T],H^{s+9}(R)),f\in L^2([0,T],H^{s})$.\\

Now we prove (\ref{3.14}) for 
$$u\in L^\infty([0,T],H^{s+\frac{1}{2}(R)})\cap C^1([0,T],H^{r}(R))\cap L^2([0,T], H^{s+2}(R)),$$ 
and $f\in L^2([0,T],H^{s-1})$.
Let $\phi\in C_0^\infty(R)$ and $\phi^\delta(\xi)=\frac{1}{\delta}\phi (\frac{\xi}{\delta})$, we have $u*\phi^\delta \in C^1([0,T],H^{s+9}(R))$ and $u*\phi^\delta$ satisfies the following equation
\begin{multline}\label{3.28}
(u*\phi^\delta)_t+b \partial_\xi^3 H[(u*\phi^\delta)]+\epsilon \partial_\xi^6(u*\phi^\delta)\\
=f*\phi^\delta+\left[b\partial^3_\xi H[u*\phi^\delta]-\left(bH[\partial_\xi^3 u]\right)*\phi^\delta\right].
\end{multline}
Note that $bH[\partial_\xi^3 u],\partial^{s-1}_\xi \left(bH[\partial_\xi^3 u]\right)\in L^2([0,T],R)$, we have
\begin{equation}\label{3.29}
\int_0^T \norm{\left(b\partial^{3}_\xi H[u]\right)*\phi^\delta-b\partial^{3}_\xi H[u]}^2_{H^{s-1}}dt\to 0, \text{ as } \delta\to 0;
\end{equation}
and
\begin{multline}\label{3.30}
\int_0^T \norm{b\partial^3_\xi H[u*\phi^\delta]-b\partial^3_\xi H[u]}^2_{H^{s-1}}dt\\
\leq 2\norm{b}_{L^\infty([0,T], H^s(R))}\int_0^T\norm{u*\phi^\delta-u}^2_{H^{s+2}}dt;
\end{multline}
so
\begin{equation}\label{3.31}
\int_0^T \norm{b\partial^3_\xi H[u*\phi^\delta]-b\partial^3_\xi H[u]}_{H^{s-1}}dt\to 0, \text{ as } \delta\to 0;
\end{equation}
therefore
\begin{equation}\label{3.32}
\int_0^T \norm{b\partial^3_\xi H[u*\phi^\delta]-(bH[\partial_\xi^3 u])*\phi^\delta}_{H^{s-1}}dt\to 0, \text{ as } \delta\to 0.
\end{equation}

Now let 
\begin{equation*}
f^\delta =f*\phi^\delta+[b\partial^3_\xi H[u*\phi^\delta]-(bH[\partial_\xi^3 u])*\phi^\delta]
\end{equation*}
 applying (\ref{3.14}) to $u*\phi^\delta$, 
we have
\begin{multline}\label{3.33}
\norm{u*\phi^\delta(\cdot,t)}^2_{H^{s-j+1/2}(R)}+\beta\int_0^t e^{K(t-\rho)}\norm{\partial_\xi^{s-j+2}u*\phi^\delta(\cdot,\rho)}^2_{L^2(R)}d\rho\\
+\epsilon \int_0^t e^{K(t-\rho)}\norm{u*\phi^\delta}^2_{H^{s-j+3\frac{1}{2}}(R)}d\rho \\
\leq e^{Kt}\norm{u_0*\phi^\delta(\cdot,t)}^2_{H^{s-j+1/2}(R)}+K_1\int_0^t e^{K(t-\rho)}\norm{f^\delta(\cdot,\rho)}^2_{H^{s-j-1}(R)}d\rho,
\end{multline}
Let $\delta\to 0$, we have (\ref{3.14}).\\
We now prove that $u\in C([0,T],H^{s+\frac{1}{2}}(R))$. Applying (\ref{3.27}) with $j=0$ to $u*\phi^\delta$, we have
\begin{multline}\label{3.34}
\frac{d}{dt}\norm{u*\phi^\delta}^2_{H^{s+\frac{1}{2}}}+\epsilon\norm{u*\phi^\delta}^2_{H^{s+3\frac{1}{2}}}
+\beta\norm{\partial_\xi^{s+2}u*\phi^\delta}^2_{L^{2}}\\
\leq K_1\norm{f^\delta}_{H^{s-1}}+K\norm{u*\phi^\delta}_{H^{s+\frac{1}{2}}}
\end{multline}
Integrating both sides of (\ref{3.34}) with respect to $t$ from $t_0$ to $t$, and then let $\delta\to 0$, we obtain
\begin{equation}\label{3.35}
\left\vert \norm{u(t)}^2_{H^{s+\frac{1}{2}}}-\norm{u(t_0)}^2_{H^{s+\frac{1}{2}}}\right\vert\leq \left\vert \int_{t_0}^t g(t)dt\right\vert
\end{equation}
for $g(t)\in L^1[0,T]$. This shows that $\norm{u(t)}^2_{H^{s+\frac{1}{2}}}\to\norm{u(t_0)}^2_{H^{s+\frac{1}{2}}}$ as $t\to t_0$. Since $u(t)\to u(t_0)$ weakly in $H^{s+\frac{1}{2}}$, therefore
\begin{equation}\label{3.36}
\norm{u(t)-u(t_0)}^2_{H^{s+\frac{1}{2}}(R)}\to 0 \text{ as }t \to t_0.
\end{equation}
\end{proof}

\begin{lemma}
Let $f\in L^2([0,T], H^{s-1}(R))\cap C([0,T], H^{s-5/2}(R))$ and $u_0\in H^{s+\frac{1}{2}}$, $b$ is as in Lemma 3.3,  then the initial value problem 
\begin{equation}\label{3.37}
\begin{split}
&u_t+b H[\partial^3_\xi u]=f(t,\xi),\\
 &u|_{t=0}=u_0~~~
\end{split}
\end{equation}
has a unique solution 
$$u\in C([0,T], H^{s+\frac{1}{2}}(R))\cap C^1([0,T], H^{s-5/2}(R))\cap L^2([0,T], H^{s+2}(R)).$$
\end{lemma} 
\begin{proof}
Let $\epsilon>0$, by Lemma 3.2, there is $u^\epsilon\in C([0,T],H^{s+2}(R))\cap L^2([0,T],H^{s+5}(R))$ that is the unique solution of 
\begin{equation}\label{3.38}
\begin{split}
&u^\epsilon_t+bH[\partial^3_\xi u^\epsilon]-\epsilon\partial^6_\xi u^\epsilon=f(t,\xi),\\
 &u^\epsilon|_{t=0}=u_0*\phi^\epsilon~~~~
\end{split}
\end{equation}
From (\ref{3.38}), we have $u^\epsilon\in C^1([0,T],H^{-1}(R))$. From (\ref{3.14}) we have that $\{u^\epsilon\}$ is a bounded set in $C([0,T],H^{s+\frac{1}{2}})\cap L^2([0,T],H^{s+3\frac{1}{2}}(R))$.\\
Now
\begin{equation}\label{3.39}
\begin{split}
&(u^\epsilon-u^\delta)_t
+b(\xi)H[\partial^3_\xi (u^\epsilon-u^\delta)]-\epsilon\partial^6_\xi (u^\epsilon -u^\delta)=(\epsilon-\delta)\partial_\xi^6u^\delta,\\
 &(u^\epsilon -u^\delta)|_{t=0}=u_0*\phi^\epsilon -u_0*\phi^\delta~~~~~~~~~~~
\end{split}
\end{equation}
Applying (\ref{3.14}) with $j=0$ to (\ref{3.39}), we have
\begin{multline}\label{3.45}
\norm{u^\epsilon(\cdot,t)-u^\delta}^2_{H^{s+1/2}(R)}
+\beta \int_0^t e^{K(t-\rho)}\norm{\partial_\xi^{s+2}(u^\epsilon(\cdot,\rho)-u^\delta)}^2_{L^2(R)}d\rho\\
+\epsilon \int_0^t e^{K(t-\rho)}\norm{u^\epsilon-u^\delta}^2_{H^{s+3\frac{1}{2}}(R)}d\rho
\\
\leq e^{Kt}\norm{u_0(\cdot,t)*(\phi^\epsilon-\phi^\delta)}^2_{H^{s+1/2}(R)}+K_1(\epsilon-\delta)^2\int_0^t e^{K(t-\rho)}\norm{\partial^{s+5}u^\delta}^2_{L^2(R)}d\rho,
\end{multline}
Therefore $\{u^\epsilon\}$ is a Cauchy sequence in $C([0,T],H^{s+\frac{1}{2}})\cap L^2([0,T],H^{s+2}(R))$ and there exist a unique $u\in C([0,T],H^{s+\frac{1}{2}})\cap L^2([0,T],H^{s+2}(R))$ such that $u^\epsilon\to u $ 
in $C([0,T],H^{s+\frac{1}{2}})\cap L^2([0,T],H^{s+2}(R))$. Letting $\epsilon\to 0$ in (\ref{3.38}) , we obtain that $u$ is a  solution of (\ref{3.37}). $u\in C^1([0,T],H^{s-5/2})$ follows from (\ref{3.37}) and  $f\in C^1([0,T],H^{s-5/2})$.
\end{proof}

\section{The quasilinear equation}

Now we rewrite equation (\ref{1.40}) as a quasilinear equation. From (\ref{1.41}), we have
\begin{equation}\label{4.1}
Q[u]=\tau e^{-2(h^I+u)}H\left[Q_1[u]H[\partial^2_\xi u]\right]+Q_2[u];
\end{equation}
where
\begin{equation}\label{4.2}
Q_1[u]=\left(1-\gamma\cos( 4q^I-4H[u])\right) e^{-(h^I+u)};
\end{equation}
\begin{equation}
\label{4.3}
Q_2[u]= e^{-2(h^I+u)}\left\{ 1+\tau H\left[Q_1[u]H[\partial^2_\xi h^I]\right]+\tau H\left[(\partial_\xi Q_1[u])H[\partial_\xi (h^I+u)]\right]\right\}.
\end{equation}
Using commutator operator, we have
\begin{equation}\label{4.4}
H[Q_1[u]H[\partial^2_\xi u]]=-Q_1[u] \partial^2 u+[H,Q_1[u]]H[\partial^2_\xi u];
\end{equation}
so (\ref{4.1}) becomes
\begin{equation}\label{4.5}
Q[u]=-B[u]\partial_\xi^2 u +Q_2[u]+Q_3[u];
\end{equation}
where 
\begin{equation}\label{4.6}
B[u]=\tau e^{-2(h^I+u)}Q_1[u]=\left(1-\gamma\cos( 4q^I-4H[u])\right) e^{-3(h^I+u)};
\end{equation}
\begin{equation}\label{4.7}
Q_3[u]=\tau e^{-2(h^I+u)}[H,Q_1[u]]H[\partial^2_\xi u].
\end{equation}
We note that $Q_2[u]$ and $Q_3[u]$ consist of lower order terms compared to $\partial^2_\xi u$ since the commutator operator is smoothing operator from Lemma 2.8.\\
From (\ref{4.6}), we have
\begin{equation}\label{4.8} 
B[u]\geq \tau (1-\gamma)=\beta>0.
\end{equation}
Taking derivative in (\ref{4.5}), we have
\begin{equation}\label{4.9}
\partial_\xi Q[u]=- B[u]\partial_\xi^3 u +Q_4[u];
\end{equation}
where
\begin{equation}\label{4.10}
Q_4[u]=(\partial_\xi B[u])\partial_\xi^2 u+\partial_\xi (Q_2[u]+Q_3[u]).
\end{equation}
Using commutator operator, we have
\begin{equation}\label{4.11}
H[B[u]\partial^3u]=B[u]H[\partial_\xi^3u]+[H,B[u]]\partial_\xi^3u,
\end{equation}
and we can write
\begin{equation}\label{4.12}
H[\partial_\xi Q[u]]=-B[u]H[\partial_\xi^3u]+Q_5[u];
\end{equation}
where
\begin{equation}\label{4.13}
Q_5[u]=-[H,B[u]]\partial_\xi^3u+H[Q_4[u]].
\end{equation}
We note that $Q_5[u]$ consist of lower order terms compared to $\partial^3_\xi u$ since the commutator operator is smoothing operator.\\
Now using (\ref{4.12}), (\ref{1.40}) can be written as
\begin{equation}\label{4.14}
u_t+B[u]H[\partial^3u]=\mathcal{N}[u];
\end{equation}
where
\begin{equation}\label{4.15}
\mathcal{N}[u]=(h^I_\xi+u_\xi)H[Q[u]]-(q^I_\xi-H[u])Q[u]+Q_5[u];
\end{equation}
 and $\mathcal{N}[u]$ consist of lower order terms compared to $\partial^3_\xi u$ .

\begin{lemma}
If $u\in  C([0,T], H^{s+1/2}(R)), s\ge 3$, $\sup_{t\in [0,T]}\norm{u}_{H^{s=\frac{1}{2}}(R)}\leq M$, then
$Q_1[u]\in C([0,T], H^{s+1/2}(R)),B[u]\in C([0,T], H^{s+1/2}(R))$ and  for any $t\in [0,T]$,
\begin{equation}\label{4.16}
\norm{Q_1[u]}_{H^{s+1/2}(R)}\leq C(M),~~\norm{B[u]}_{H^{s+1/2}(R)}\leq C(M),
\end{equation}
where $C(M)$ depends on $M$, $\tau$ and $\gamma$.
\end{lemma}
\begin{proof}
Note that $e^{-h^I}\in H^\infty, q^I\in \dot{H}^\infty$ from (\ref{1.37}). Applying Lemma 2.2 and Lemma 2.4 to (\ref{4.2}) and (\ref{4.6}), we obtain the lemma.
\end{proof}
\begin{lemma}
If $u\in  C([0,T], H^{s+1/2}(R)), s\ge 4$,  $\sup_{t\in [0,T]}\norm{u}_{H^{s+\frac{1}{2}}(R)}\leq M$, then 
$Q_2[u]\in C([0,T], H^{s-1/2}(R))$, $Q_3[u]\in C([0,T], H^{s+1/2}(R))$,\\
 $Q[u]\in C([0,T], H^{s-3/2}(R)),$ and for any $t\in [0,T]$,
\begin{equation}\label{4.17}
\begin{split}
&\norm{Q_2[u]}_{H^{s-1/2}(R)}\leq C(M),~~\norm{Q_3[u]}_{H^{s+1/2}(R)}\leq C(M),\\
&\norm{Q[u]}_{H^{s-3/2}(R)}\leq C(M).
\end{split}
\end{equation}
\end{lemma}
\begin{proof}
Applying Lemma 2.2-Lemma 2.5 and Lemma 4.1 to (\ref{4.3}), we obtain $Q_2[u]\in C([0,T], H^{s-1/2}(R))$.\\
From Lemma 2.8, Lemma 2.4,
\begin{multline}\label{4.18}
\left\vert\left\vert[H,Q_1[u]]H[\partial^2_\xi u]\right\vert\right\vert_{H^{s+1/2}}\\
\leq C(M)\norm{Q_1[u]}_{H^{s+1/2}}\norm{H[\partial^2_\xi u]}_{H^{s-3/2}}\leq C(M)\norm{u}^2_{H^{s+1/2}}.
\end{multline}
The lemma follows from (\ref{4.7}) and (\ref{4.5}) and Lemma 4.1.
\end{proof}

\begin{lemma}
Let $T>0,s\ge 4, u\in C([0,T],H^{s+\frac{1}{2}}(R))$ and  $\sup_{t\in [0,T]}\norm{u}_{H^{s+\frac{1}{2}}(R)}\leq M$, then $Q_4[u]\in C([0,T], H^{s-3/2}(R)),Q_5[u]\in C([0,T], H^{s-3/2}(R))$,\\
$\mathcal{N}[u]\in C([0,T], H^{s-3/2}(R))$ and for every $t\in [0,T]$,
\begin{equation}\label{4.19}
\begin{split}
&\norm{Q_4[u]}_{H^{s-3/2}(R)}\leq C(M),~~\norm{Q_5[u]}_{H^{s-3/2}(R)}\leq C(M),~\\
&\norm{\mathcal{N}[u]}_{H^{s-3/2}(R)}\leq C.
\end{split}
\end{equation}
\end{lemma}
\begin{proof}
$Q_4[u]\in C([0,T], H^{s-3/2}(R))$ follows from (\ref{4.10}) and lemma 4.1 and Lemma 4.2.\\
From Lemma 2.8, Lemma 4.1,
\begin{equation}\label{4.20}
\left\vert\left\vert[H,B[u]]\partial^3_\xi u\right\vert\right\vert_{H^{s-1/2}}\leq C \norm{B[u]}_{H^{s-1/2}}\norm{\partial^3_\xi u}_{H^{s-5/2}}\leq C(M) \norm{u}^2_{H^{s+1/2}}.
\end{equation}
 The above inequality and (\ref{4.13}) imply that $Q_5[u]\in C([0,T], H^{s-3/2}(R)) $,  and $\mathcal{N}[u]\in C([0,T], H^{s-3/2}(R))$ follows from (\ref{4.15}) and Lemma 4.2.
\end{proof}

\begin{lemma}
If $u\in L^2([0,T], H^{s+2}(R))\cap C([0,T], H^{s+\frac{1}{2}}(R)), s\ge 4$, and\\
 $\sup_{t\in [0,T]}\norm{u}_{H^{s+\frac{1}{2}}(R)}\leq M$, then for a.e $t\in [0,T]$,
\begin{multline}\label{4.21}
\norm{\partial^{s-1}_\xi Q_4[u] (\cdot,t)}^2_{L^2(R)}
\leq C(M)\int (\partial^{s+1}_\xi u(\xi, t))^2d\xi +C(M),
\end{multline}
where $C(M)$ is a constant depending on M.
\end{lemma}
\begin{proof}
By (\ref{4.10}), we have
\begin{equation}\label{4.22}
\partial_\xi^{s-1}Q_4[u]=\partial_\xi^{s-1}\left[(\partial_\xi B[u])\partial_\xi^2 u\right]+\partial^s_\xi [Q_2[u]+Q_3[u]].
\end{equation}
We can write
\begin{equation}\label{4.23}
\partial_\xi^{s-1}\left[(\partial_\xi B[u])\partial_\xi^2 u\right]=(\partial_\xi B[u])\partial_\xi^{s+1}u+\sum_{k=1}^{s-1}(\partial_\xi^{k+1}B[u])\partial_\xi^{s+1-k}u.
\end{equation}
The first term in (\ref{4.23}) can be estimated as
\begin{equation}\label{4.24}
\int[(\partial_\xi B[u])\partial_\xi^{s+1}u]^2d\xi\leq \norm{(\partial_\xi B[u])}^2_{L^\infty}\int[\partial_\xi^{s+1}u]^2d\xi\leq C(M)\int[\partial_\xi^{s+1}u]^2d\xi;
\end{equation}
similarly, the terms in the sum of (\ref{4.23}) can be bounded by C(M) from Lemma 4.1.\\
By examining (\ref{4.3}), we see that the highest order term in $\partial^sQ_2[u]$ comes from $H[\partial_\xi^{s+1}Q_1[u]]$, which is, from (\ref{4.2}), $H\left[Q_1[u]\partial_\xi^{s+1}uH[\partial_\xi(h^I+u)]\right]$ that can be estimated as
\begin{multline}\label{4.25}
\int\left( H\left[Q_1[u]\partial_\xi^{s+1}uH[\partial_\xi(h^I+u)]\right]\right)^2d\xi=\int\left( Q_1[u]\partial_\xi^{s+1}uH[\partial_\xi(h^I+u)]\right)^2d\xi\\
\leq \norm{Q_1[u]H[\partial_\xi(h^I+u)]}^2_{L^\infty}\int[\partial_\xi^{s+1}u]^2d\xi\leq C(M)\int[\partial_\xi^{s+1}u]^2d\xi;
\end{multline}
the lower order terms in $\partial^sQ_2[u]$ can be estimated by $C(M)$ from Lemma 4.1.\\
Now from (\ref{4.7})
\begin{equation}\label{4.26}
\partial_\xi^sQ_3[u]=\tau\sum_{k=0}^{s} \partial_\xi^k\left(e^{-2(h^I+u)}\right)\partial_\xi^{s-k}\left([H,Q_1[u]]H[\partial^2_\xi u]\right).
\end{equation}
so
\begin{multline}\label{4.27}
\int\left(\partial_\xi^sQ_3[u]\right)^2d\xi\leq 2\tau\sum_{k=0}^{s}\int \left(\partial_\xi^ke^{-2(h^I+u)}\right)^2\left[\partial_\xi^{s-k}\left([H,Q_1[u]]H[\partial^2_\xi u]\right)\right]^2d\xi\\
\leq 2\tau\sum_{k=0}^{s}\left\{\int \left(\partial_\xi^ke^{-2(h^I+u)}\right)^4d\xi\right\}^{1/2}\left\{\int\left[\partial_\xi^{s-k}\left([H,Q_1[u]]H[\partial^2_\xi u]\right)\right]^4d\xi\right\}^{1/2}\\
\leq 2K^2_0\tau\sum_{k=0}^{s}\norm{e^{-2(h^I+u)}}^2_{H^{k+1/4}}\norm{[H,Q_1[u]]H[\partial^2_\xi u]}^2_{H^{s-k+1/4}}\\
\leq C\sum_{k=0}^{s}\norm{e^{-2(h^I+u)}}^2_{H^{k+1/4}}\norm{Q_1[u]}^2_{H^{s-k+1/4}}\norm{H[\partial^2_\xi u]}^2_{H^{s-k+1/4-2}}\\
\leq C \sum_{k=0}^{s}\norm{u}^2_{H^{k+1/4}}\norm{u}^2_{H^{s-k+1/4}}\norm{u}^2_{H^{s-k+1/4}}\leq C(M).~~~~~~~~~~~~~~~~~~~~~~~~~~~~~~~~~~~~~~~~~~~~~~~~~~~~~~~~~~~~~~~~~
\end{multline}
The lemma follows from (\ref{4.22})-(\ref{4.27}).
\end{proof}
\begin{lemma}
If $u\in L^2([0,T], H^{s+2}(R))\cap C([0,T], H^{s+\frac{1}{2}}(R)),s\ge 4$, and\\
$\sup_{t\in [0,T]}\norm{u}_{H^{s+\frac{1}{2}}(R)}\leq M$, then for a.e $t\in [0,T]$,
\begin{multline}\label{4.28}
\norm{\partial^{s-1}_\xi Q_5[u] (\cdot,t)}^2_{L^2(R)}
\leq C(M)\int (\partial^{s+1}_\xi u(\xi, t))^2d\xi +C(M),
\end{multline}
where $C(M)$ is a constant depending on M.
\end{lemma}
\begin{proof}
from (\ref{4.13}), we have
\begin{equation}\label{4.29}
\partial_\xi^{s-1}Q_5[u]=-\partial_\xi^{s-1}[H,B[u]]\partial_\xi^3u+\partial_\xi^{s-1}H[Q_4[u]];
\end{equation}
the first term on the right side of (\ref{4.29}) can be estimated as
\begin{multline}\label{4.30}
\int\left(\partial_\xi^{s-1}[H,B[u]]\partial_\xi^3u\right)^2d\xi\leq \norm{[H,B[u]]\partial_\xi^3u}_{H^{s-1}}\\
\leq C\norm{B[u]}_{H^{s-1}}\norm{\partial_\xi^3u}_{H^{s-3}}\leq C(M)
\end{multline}
Lemma 2.8 and Lemma 4.1 were used in obtaining (\ref{4.30}).\\

The lemma follows from (\ref{4.29}) and Lemma 4.4.
\end{proof}

\begin{lemma}
If $u\in L^2([0,T], H^{s+2}(R))\cap C([0,T], H^{s+\frac{1}{2}}(R)),s\ge 4$, and\\
 $\sup_{t\in [0,T]}\norm{u}_{H^{s+\frac{1}{2}}(R)}\leq M$, then for a.e $t\in [0,T]$,
\begin{equation}\label{4.31}
\norm{\partial^{s-1}_\xi N[u] (\cdot,t)}^2_{L^2(R)}
\leq C(M)\int (\partial^{s+1}_\xi u(\xi, t))^2d\xi +C(M),
\end{equation}
where $C(M)$ is a constant depending on M.
\end{lemma}
\begin{proof}
From (\ref{4.15}) and (\ref{4.5}), we have
\begin{multline}\label{4.32}
\mathcal{N}[u]=-(h^I_\xi+u_\xi)H[B[u]\partial_\xi^2 u]+(h^I_\xi+u_\xi)H[Q_2[u]+Q_3[u]]\\
+(q^I_\xi-H[u])B[u]\partial_\xi^2 u-(q^I_\xi-H[u])(Q_2[u]+Q_3[u])+Q_5[u];
\end{multline}
Note that two highest order terms in $\partial_\xi^{s-1}\mathcal{N}[u]$ can be estimated in the same fashion as
\begin{multline}\label{4.33}
\int\left( (h^I_\xi+u_\xi)H[B[u]\partial_\xi^{s+1} u]\right)^2d\xi\leq \norm{(h^I_\xi+u_\xi)}^2_{L^\infty}\int\left(H[B[u]\partial_\xi^{s+1} u]\right)^2d\xi\\
\leq \norm{(h^I_\xi+u_\xi)}^2_{L^\infty}\int\left(B[u]\partial_\xi^{s+1} u\right)^2d\xi\\
\leq \norm{(h^I_\xi+u_\xi)}^2_{L^\infty}\norm{B[u]}^2_{L^\infty}\int\left(\partial_\xi^{s+1} u\right)^2d\xi\leq C(M)\int\left(\partial_\xi^{s+1} u\right)^2d\xi.\\
\end{multline}
The lemma then follows from (\ref{4.32}) and Lemma 4.5.

\end{proof}

\begin{lemma}
Assume that $u^k\in L^\infty([0,T], H^{s+\frac{1}{2}}(R)), s\ge 4, k=1,2$; and \\
$\sup_{t\in [0,T]}\norm{u^k(\cdot,t)}_{H^{s+\frac{1}{2}}(R)}\leq M$, then for any $t\in [0,T]$ and $\frac{1}{2}<r\leq s+\frac{1}{2}$,
\begin{equation}\label{4.34}
\norm{Q_1[u^1](\cdot,t)-Q_1[u^2](\cdot,t)}^2_{H^{r}(R)}\leq C\norm{u^1(\cdot,t)-u^2(\cdot,t)}^2_{H^{r}(R)},
\end{equation}
\begin{equation}\label{4.35}
\norm{B[u^1](\cdot,t)-B[u^2](\cdot,t)}^2_{H^{r}(R)}\leq C\norm{u^1(\cdot,t)-u^2(\cdot,t)}^2_{H^{r}(R)},
\end{equation}
\end{lemma}
\begin{proof}
From (\ref{4.2}), we have
\begin{multline}\label{4.36}
Q_1[u^1]-Q_1[u^2]=
\left(1-\gamma\cos (4q^I-4H[u_1])\right)e^{-2u^2-2h^I}(e^{2(u^1-u^2)}-1)\\
+2\gamma e^{-2u^1-2h^I}\sin (2H[u^1+u^2])\sin (2H[u^1-u^2]).
\end{multline}
Using Lemma 2.1 and Lemma 2.2 to obtain
\begin{multline}\label{4.37}
\norm{Q_1[u^1]-Q_1[u^2]}_{H^r}\leq
\norm{\left(1-\gamma\cos (4q^I-4H[u_1])\right)e^{-2u^2-2h^I}(e^{2(u^1-u^2)}-1)}_{H^r}\\
+2|\gamma| \norm{ e^{-2u^1-2h^I}\sin (2H[u^1+u^2])\sin (2H[u^1-u^2])}_{H^r}\\
\leq \norm{\left(1-\gamma\cos (4q^I-4H[u_1])\right)}_{H^r}\norm{e^{-2u^2-2h^I}}_{H^r}\norm{(e^{2(u^1-u^2)}-1)}_{H^r}\\
+2|\gamma| \norm{ e^{-2u^1-2h^I}}_{H^r}\norm{\sin (2H[u^1+u^2])}_{H^r}\norm{\sin (2H[u^1-u^2])}_{H^r}\\
\leq C(M)\norm{u^1-u^2}_{H^r}~~~~~~~~~~~~~~~~~~~~~~~~~~~
\end{multline}

(\ref{4.35}) can be proved similarly.
\end{proof}

\begin{lemma}
Assume that $u^k\in L^\infty([0,T], H^{s+\frac{1}{2}}(R))\cap L^2([0,T], H^{s+2}(R)), k=1,2, s\ge 4; \text{ and }\sup_{t\in [0,T]}\norm{u^k(\cdot,t)}_{H^{s+\frac{1}{2}}(R)}\leq M$, then for a.e $t\in [0,T]$,
\begin{multline}\label{4.38}
\norm{\partial^{s+1}_\xi Q_1[u^1](\cdot,t)-\partial^{s+1}_\xi Q_1[u^2](\cdot,t)}^2_{L^2(R)}\\
\leq C(M)\norm{u^1(\cdot,t)-u^2(\cdot,t)}^2_{H^{s+1}(R)}+C(M)\norm{u^1-u^2}_{H^{3/2}}\int(\partial_\xi^{s+1}u^2)^2d\xi,
\end{multline}
and
\begin{multline}\label{4.39}
\norm{\partial^{s}_\xi Q_2[u^1](\cdot,t)-\partial^s_\xi Q_2[u^2](\cdot,t)}^2_{L^2(R)}\\
\leq C(M)\norm{u^1(\cdot,t)-u^2(\cdot,t)}^2_{H^{s+1}(R)}+C(M)\norm{u^1-u^2}_{H^{3/2}}\int(\partial_\xi^{s+1}u^2)^2d\xi.
\end{multline}
\end{lemma}
\begin{proof}
From (\ref{4.36}), the highest order terms in $\partial^{s+1}_\xi Q_1[u^1](\cdot,t)-\partial^{s+1}_\xi Q_1[u^2]$ are  such terms as $\left(1-\gamma\cos (4q^I-4H[u_1])\right)e^{-2u^2-2h^I}(e^{2(u^1-u^2)}-1)\partial^{s+1}(u^1-u^2)$ and $\left(1-\gamma\cos (4q^I-4H[u_1])\right)e^{-2u^2-2h^I}(e^{2(u^1-u^2)}-1)\partial^{s+1}u^2$, which can be estimated as
\begin{multline}\label{4.40}
\int\left( (1-\gamma\cos (4q^I-4H[u_1]))e^{-2u^2-2h^I}(e^{2(u^1-u^2)}-1)\partial^{s+1}(u^1-u^2)\right)^2d\xi\\
\leq \norm{\left(1-\gamma\cos (4q^I-4H[u_1])\right)e^{-2u^2-2h^I}}^2_{L^\infty}\\
\times\norm{(e^{2(u^1-u^2)}-1)}^2_{L^\infty}\int\left(\partial^{s+1}_\xi(u^1-u^2)\right)^2d\xi\\
\leq C(M)\int\left(\partial^{s+1}_\xi(u^1-u^2)\right)^2d\xi;
\end{multline}
\begin{multline}\label{4.41}
\int\left( (1-\gamma\cos (4q^I-4H[u_1]))e^{-2u^2-2h^I}(e^{2(u^1-u^2)}-1)\partial^{s+1}u^2\right)^2d\xi\\
\leq \norm{\left(1-\gamma\cos (4q^I-4H[u_1])\right)e^{-2u^2-2h^I}}^2_{L^\infty}\norm{(e^{2(u^1-u^2)}-1)}^2_{L^\infty}\int\left(\partial^{s+1}_\xi u^2\right)^2d\xi\\
\leq C(M)\norm{u^1-u^2}_{H^{3/2}}\int\left(\partial^{s+1}_\xi u^2\right)^2d\xi;
\end{multline}
the lower terms in $\partial^{s+1}_\xi Q_1[u^1](\cdot,t)-\partial^{s+1}_\xi Q_1[u^2]$ can be bounded by \\
$C(M)\int\left(\partial^{s+1}_\xi (u^1-u^2)\right)^2d\xi$ from Lemma 4.7. So (\ref{4.38}) is proved.\\
From (\ref{4.3}), we obtain
\begin{multline}\label{4.42}
Q_2[u^1]-Q_2[u^2]=e^{-2(h^I+u^1)}\left(1-e^{-2(u^2-u^1)}\right)\left(1 +\tau H[Q_1[u^1]]H[\partial_\xi^2h^I]\right)\\
+e^{-2(h^I+u^2)}\left(\tau H[Q_1[u^1]-Q_1[u^2]]H[\partial_\xi^2h^I]\right)\\
+(\tau H[\partial_\xi (Q_1[u^1]-Q_1[u^2])]H[\partial_\xi (h^I+u^1)]+\tau H[\partial_\xi Q_1[u^2]]H[\partial_\xi (u^1-u^2)];
\end{multline}
(\ref{4.39}) then follows from (\ref{4.42}) and (\ref{4.38}).

\end{proof}

\begin{lemma}
Assume that $u^k\in L^\infty([0,T], H^{s+\frac{1}{2}}(R))\cap L^2([0,T], H^{s+2}(R)), k=1,2, s\ge 4; \text{ and }\sup_{t\in [0,T]}\norm{u^k(\cdot,t)}_{H^{s+\frac{1}{2}}(R)}\leq M$, then for a.e $t\in [0,T]$,
\begin{equation}\label{4.43}
\norm{\partial^{s}_\xi Q_3[u^1](\cdot,t)-\partial^s_\xi Q_3[u^2](\cdot,t)}^2_{L^2(R)}\leq C(M)\norm{u^1(\cdot,t)-u^2(\cdot,t)}^2_{H^{s}(R)},
\end{equation}
\begin{multline}\label{4.44}
\norm{\partial^{s-1}_\xi Q_4[u^1](\cdot,t)-\partial^{s-1}_\xi Q_4[u^2](\cdot,t)}^2_{L^2(R)}\\
\leq C(M)\norm{u^1(\cdot,t)-u^2(\cdot,t)}^2_{H^{s+1}(R)}+C(M)\norm{u^1-u^2}_{H^{3/2}}\int(\partial_\xi^{s+1}u^2)^2d\xi.
\end{multline}
\end{lemma}
\begin{proof}
From (\ref{4.7})
\begin{multline}\label{4.45}
Q_3[u^1]-Q_3[u^2]=\tau e^{-2(h^I+u^1)}\left(1-e^{-2(u^2-u^1)}\right)[H,Q_1[u^1]]H[\partial_\xi^2u^1]\\
+\tau e^{-2(h^I+u^2)}[H,Q_1[u^1]-Q_1[u^2]]H[\partial_\xi^2u^1]\\
+\tau e^{-2(h^I+u^2)}[H,Q_1[u^1]]H[\partial_\xi^2(u^1-u_2)].
\end{multline}
Using Lemma 2.1, Lemma 2.4 and Lemma 2.8 to obtain
\begin{multline}\label{4.46}
\norm{Q_3[u^1]-Q_3[u^2]}_{H^s}\\
\leq \tau \norm{e^{-2(h^I+u^1)}}_{H^s}\norm{\left(1-e^{-2(u^2-u^1)}\right)}_{H^s}\norm{[H,Q_1[u^1]]H[\partial_\xi^2u^1]}_{H^s}\\
+\tau \norm{e^{-2(h^I+u^2)}}_{H^s}\norm{[H,Q_1[u^1]-Q_1[u^2]]H[\partial_\xi^2u^1]}_{H^s}\\
+\tau \norm{e^{-2(h^I+u^2)}}_{H^s}\norm{[H,Q_1[u^1]]H[\partial_\xi^2(u^1-u^2)]}_{H^s}\\
\leq \tau C(M)\norm{u^2-u^1}_{H^s}\norm{Q_1[u^1]]}_{H^s}\norm{H[\partial_\xi^2u^1]}_{H^{s-2}}\\
+\tau C(M)\norm{Q_1[u^1]-Q_1[u^2]}_{H^s}\norm{H[\partial_\xi^2u^1]}_{H^{s-2}}\\
+\tau C(M)\norm{Q_1[u^1]}_{H^s}\norm{H[\partial_\xi^2(u^1-u^2)]}_{H^{s-2}}\\
\leq C(M)\norm{(u^1-u^2)]}_{H^{s}}.~  
\end{multline}
From (\ref{4.10}), we have
\begin{multline}\label{4.47}
Q_4[u^1]-Q_4[u^2]=(\partial_\xi B[u^1])(\partial_\xi^2u^1-\partial_\xi^2u^2)+(\partial_\xi B[u^1]-\partial_\xi B[u^2])\partial_\xi^2u^2\\
+\partial_\xi \left(Q_2[u^1]-Q_2[u^2]\right)+\partial_\xi \left(Q_3[u^1]-Q_3[u^2]\right).
\end{multline}
Note that the highest order term from the first two terms on the right hand side of (\ref{4.47}) can
be estimated as
\begin{multline}\label{4.48}
\int\left((\partial_\xi B[u^1])(\partial_\xi^{s+1}u^1-\partial_\xi^{s+1}u^2)\right)^2d\xi\\
\leq \norm{(\partial_\xi B[u^1])}^2_{L^\infty}\int\left(\partial_\xi^{s+1}(u^1-u^2)\right)^2d\xi\\
\leq C(M)\int\left(\partial_\xi^{s+1}(u^1-u^2)\right)^2d\xi
\end{multline}
and
\begin{multline}\label{4.49}
\int\left(\partial_\xi B[u^1]-\partial_\xi B[u^2])\partial_\xi^{s+1}u^2\right)^2d\xi\\
\leq \norm{(\partial_\xi B[u^1]-\partial_\xi B[u^2])}^2_{L^\infty}\int\left(\partial_\xi^{s+1}u^2\right)^2d\xi\\
\leq C(M)\norm{u^1-u^2}_{H^{3/2}}\int\left(\partial_\xi^{s+1}u^2\right)^2d\xi.
\end{multline}
Then (\ref{4.44}) follows from (\ref{4.47})-(\ref{4.49}), ( \ref{4.43}) and Lemma 4.8.
\end{proof}
\begin{lemma}
Assume that $u^k\in L^\infty([0,T], H^{s+\frac{1}{2}}(R))\cap L^2([0,T], H^{s+2}(R)), k=1,2, s\ge 4; \text{ and }\sup_{t\in [0,T]}\norm{u^k(\cdot,t)}_{H^{s+\frac{1}{2}}(R)}\leq M$, then for a.e $t\in [0,T]$,
\begin{multline}\label{4.50}
\norm{\partial^{s-1}_\xi Q_5[u^1](\cdot,t)-\partial^{s-1}_\xi Q_5[u^2](\cdot,t)}^2_{L^2(R)}\\
\leq C(M)\norm{u^1(\cdot,t)-u^2(\cdot,t)}^2_{H^{s+1}(R)}+C(M)\norm{u^1-u^2}_{H^{3/2}}\int(\partial_\xi^{s+1}u^2)^2d\xi.
\end{multline}
\end{lemma}
\begin{proof}
From (\ref{4.13}), we have
\begin{multline}\label{4.51}
Q_5[u^1]-Q_5[u^2]=H[Q_4[u^1]-Q_4[u^2]]\\
-[H,B[u^1]-B[u^2]]\partial_\xi^3u^1-[H,B[u^2]]\partial_\xi^3(u^1-u^2).
\end{multline}
Using Lemma 2.8 and Lemma 4.7 to obtain
\begin{multline}\label{4.52}
\norm{[H,B[u^1]-B[u^2]]\partial_\xi^3u^1}_{H^{s-1}}\leq C\norm{B[u^1]-B[u^2]]}_{H^{s-1}}\norm{\partial_\xi^3u^1}_{H^{s-3}}\\
\leq C(M)\norm{u^1-u^2}_{H^{s-1}}
\end{multline}
and 
\begin{multline}\label{4.53}
\norm{[H,B[u^2]]\partial_\xi^3(u^1-u^2)}_{H^{s-1}}\leq C\norm{B[u^2]]}_{H^{s-1}}\norm{\partial_\xi^3(u^1-u^2}_{H^{s-3}}\\
\leq C(M)\norm{u^1-u^2}_{H^{s}}
\end{multline}
then (\ref{4.50}) follows from (\ref{4.51})-(\ref{4.53}) and Lemma 4.9.
\end{proof}
\begin{lemma}
Assume that $u^k\in L^\infty([0,T], H^{s+\frac{1}{2}}(R))\cap L^2([0,T], H^{s+2}(R)), k=1,2; s\ge 4, \text{ and }\sup_{t\in [0,T]}\norm{u^k(\cdot,t)}_{H^{s+\frac{1}{2}}(R)}\leq M$, then for a.e $t\in [0,T]$,
\begin{multline}\label{4.54}
\norm{\mathcal{N}[u^1](\cdot,t)-\mathcal{N}[u^2](\cdot,t)}^2_{H^{s-1}(R)}\\
\leq C(M)\norm{u^1(\cdot,t)-u^2(\cdot,t)}^2_{H^{s+1}(R)}+C(M)\norm{u^1-u^2}_{H^{3/2}}\int(\partial_\xi^{s+1}u^2)^2d\xi.
\end{multline}
\end{lemma}
\begin{proof}
From (\ref{4.15}), we obtain
\begin{multline}\label{4.55}
\mathcal{N}[u^1]-\mathcal{N}[u^2]=-(h^I_\xi+u^1_\xi)H\left[(B[u^1]-B[u^2])\partial_\xi^2u^2\right]\\
-(h^I_\xi+u^1_\xi)H\left[B[u^2]\partial_\xi^2(u^1-u^2)\right]-(u^1_\xi-u^2_\xi)H\left[B[u^2]\partial_\xi^2u^2\right]\\
+(h^I_\xi+u^1_\xi)H\left[(Q_2[u^1]-Q_2[u^2])+(Q_3[u^1]-Q_3[u^2])\right]
+(u^1_\xi-u^2_\xi)H\left[Q_2[u^2])+Q_3[u^2]\right]\\
+(q_\xi^I-H[u^1_\xi])(B[u^1]-B[u^2])\partial_\xi^2u^1+(q_\xi^I-H[u^1_\xi])B[u^2]\partial_\xi^2(u^1-u^2)\\
-(q_\xi^I-H[u^1_\xi])\left[(Q_2[u^1]-Q_2[u^2])+(Q_3[u^1]-Q_3[u^2])\right]\\
-(u^1_\xi-u^2_\xi)\left[Q_2[u^2])+Q_3[u^2]\right]+(Q_5[u^1]-Q_5[u^2]);
\end{multline}
then the lemma follows from Lemma 4.7-4.10.
\end{proof}
Now we are ready to prove the main theorem.\\

{\bf Proof of Theorem 2.1 }: We construct a sequence of functions $\{u^k\}$ by solving
\begin{equation}\label{4.56}
\begin{split}
&\partial_t u^{k+1}+B[u^k]H(\partial_\xi^3 u^{k+1})=\mathcal{N}[u^k], \\
&u^{k+1}|_{t=0}=u_0, k=0,1,2 \cdots~~
\end{split}
\end{equation}
and $u^0=0$.

 From Lemma 3.4, we have $u^k\in C([0,T], H^{s+\frac{1}{2}}(R))\cap C^1([0,T], H^{s-5/2}(R))\cap L^2([0,T], H^{s+2}(R))$ for any $T>0$ and $k=1,2 \cdot \cdot \cdot$.\\
  We are going to show 

\begin{lemma} There exists $T=T(M_0)>0$ so that for all $k$ ,
\begin{equation}\label{4.57}
\begin{split}
&\sup_{t\in [0,T]}\norm{u^k(\cdot, t)}^2_{H^{s+\frac{1}{2}}(R)}\leq 8M_0^2\\
&|\beta|\int_0^T\int |\partial_\xi^{s+2}u^k(\xi,t)|^2d\xi dt\leq 8M_0^2
\end{split}
\end{equation}
\end{lemma}
\begin{proof}
We assume that (\ref{4.57}) hold for $u^k$. From Lemma 3.3, we have
\begin{multline}\label{4.58}
\norm{u^{k+1}(\cdot,t)}^2_{H^{s+\frac{1}{2}}(R)} +|\beta|\int_0^t e^{K(t-\rho)}\norm{\partial^{s+2}_\xi u^{k+1}(\cdot,\rho)}^2_{L^2(R)}d\rho\\
\leq e^{K t}\norm{u_0}^2_{H^{s+\frac{1}{2}}(R)}+K_1\int_0^t e^{K(t-\rho)}\norm{\mathcal{N}[u^k]}^2_{H^{s-1}(R)} d\rho \\
\leq e^{K T} M_0^2+K_1e^{K T}\int_0^t \norm{\mathcal{N}[u^k]}^2_{H^{s-1}(R)} d\rho.
\end{multline}
From Lemma 4.6, we have
\begin{multline}\label{4.59}
\int_0^T \norm{\mathcal{N}[u^k]}^2_{H^{s-1}(R)} dt
\leq C(M_0)\int_0^T\int |\partial^{s+1}_\xi u^k|^2 d\xi dt +TC(M_0)
\end{multline}
From Lemma 2.7, we have for any $\eta>0$
\begin{equation}\label{4.60}
\int |\partial^{s+1}_\xi u^k|^2 d\xi dt\leq \eta \int |\partial^{s+2}_\xi u^k|^2 d\xi+K(\eta)\norm{u^k}_{H^{s+1/2}}
\end{equation}
by (\ref{4.57})- (\ref{4.60}) imply
\begin{multline}\label{4.61}
\norm{u^{k+1}(\cdot,t)}^2_{H^{s+\frac{1}{2}}(R)} +|\beta|\int_0^t e^{K(t-\rho)}\norm{\partial^{s+2}_\xi u^{k+1}(\cdot,\rho)}^2_{L^2(R)}d\rho\\
\leq e^{K T} M_0^2+8K_1e^{K T}C(M_0)\eta +8K(\eta)K_1e^{KT}C(M_0)TM_0^2+Te^{KT}K_1C(M_0)
\end{multline}
Now choose $\eta=\frac{1}{ 8K_1e^{K T}C(M_0)}$ in (\ref{4.61})  we have 
\begin{multline}\label{4.62}
\sup_{t} \norm{u^{k+1}(\cdot,t)}^2_{H^{s+\frac{1}{2}}(R)} +\beta\int_0^T \norm{\partial^{s+2}_\xi u^{k+1}(\cdot,t)}^2_{L^2(R)}dt\\
\leq e^{KT}M_0^2 +M_0^2 +8K(\eta)K_1e^{K T}TC(M_0)M_0^2 +C(M_0)K_1e^{K T} T.
\end{multline}
Then choose T so that $$e^{\beta T}\leq 3,~~8C(M_0)K(\eta)K_1e^{2K T}T\leq 2,~~C(M_0)K_1e^{KT} T\leq 2M_0^2,$$
so we have
\begin{equation*}
\sup_{t} \norm{u^{k+1}(\cdot,t)}^2_{H^{s+/2}(R)} +|\beta|\int_0^t \norm{\partial^{s+2}_\xi u^{k+1}(\cdot,\rho)}^2_{L^2(R)}d\rho\leq 8M_0^2.
\end{equation*}
Then the Lemma follows from induction.
\end{proof}

Now we consider the convergence of $\{u^k\}$. From (\ref{4.56}), we obtain
\begin{equation}\label{4.63}
\begin{split}
&\partial_t (u^{k+1}-u^k)+B[u^k] H(\partial_\xi^3 (u^{k+1}-u^k))\\
&~~=(\mathcal{N}[u^{k+1}]-\mathcal{N}[u^k]) - (B[u^k]-B[u^{k-1}])H(\partial_\xi^3u^k); \\
&(u^{k+1}-u^k)|_{t=0}=0.
\end{split}
\end{equation}
Applying Lemma 3.3 with $j=1$, we obtain 
\begin{multline}\label{4.64}
\norm{u^{k+1}-u^k}^2_{H^{s-\frac{1}{2}}(R)}+|\beta|\int_0^t \norm{\partial^{s+1}_\xi (u^{k+1}-u^k)}^2_{L^2(R)} d\rho\\
\leq K_1e^{K T}\{ \int_0^t\norm{\mathcal{N}[u^{k}]-\mathcal{N}[u^{k-1}]}^2_{H^{s-2}} d\rho\\
+\int_0^t\norm{(B[u^k]-B[u^{k-1}])H[\partial_\xi^3u^k]}^2_{H^{s-2}} d\rho\}
\end{multline}
Note that
\begin{multline}\label{4.65}
\int_0^t\norm{(B[u^k]-B[u^{k-1}])H(\partial_\xi^3u^k)}^2_{H^{s-2}} d\rho\\
\leq \int_0^t\norm{(B[u^k]-B[u^{k-1}])}^2_{H^{s-2}}\norm{\partial_\xi^3u^k}^2_{H^{s-2}} d\rho\\
\leq \sup_{t\in [0,T]} \norm{(B[u^k]-B[u^{k-1}])}^2_{H^{s-2}(R)}\int_0^T\norm{u^k}^2_{H^{s+1}(R)}\\
\leq C(M_0)\sup_{t\in [0,T]} \norm{(u^k-u^{k-1})}^2_{H^{s-1}(R)}\int_0^T\norm{u^k}^2_{H^{s+1}(R)};
\end{multline}
Using Lemma 2.8 and Lemma 4.12, we have
\begin{multline}\label{4.66}
\int_0^T\norm{u^k}^2_{H^{s+1}(R)}\leq \eta \int_0^T\norm{u^k}^2_{H^{s+2}(R)}d\rho +K(\eta)\int_0^T\norm{u^k}^2_{L^2}d\rho\\
\leq 8M_0^2\eta +8M_0^2K(\eta)T;~~~~~~~~~~~~~~~~~
\end{multline}
so
\begin{multline}\label{4.67}
\int_0^t\norm{(B[u^k]-B[u^{k-1}])H(\partial_\xi^3u^k)}^2_{H^{s-2}} d\rho\\
\leq C(M_0)\left(8M_0^2\eta +8M_0^2K(\eta)T\right)\sup_{t\in [0,T]} \norm{(u^k-u^{k-1})}^2_{H^{s-1}(R)}.
\end{multline}
Using Lemma 4.11 and Lemma 2.8 to obtain
\begin{multline}\label{4.68}
\norm{\mathcal{N}[u^{k}]-\mathcal{N}[u^{k-1}]}^2_{H^{s-2}}\\
\leq C(M_0)\norm{(u^k-u^{k-1})}^2_{H^{s}(R)}+C(M_0)\norm{(u^k-u^{k-1})}^2_{H^{3/2}(R)}\int(\partial_\xi^su^{k-1})^2d\xi
\end{multline}
Note that by Lemma 4.12
\begin{multline}\label{4.69}
\int_0^T\norm{u^{k-1}}^2_{H^{s}(R)}\leq T\sup_{t\in [0,T]}\norm{u^{k-1}}^2_{H^{s+1/2}(R)}\leq 8M_0^2T;~~~~~
\end{multline}
and Using Lemma 2.8 we have
\begin{multline}\label{4.70}
\int_0^T\norm{(u^k-u^{k-1})}^2_{H^{s}(R)}d\rho\\
\leq \eta_0\int_0^T\norm{(u^k-u^{k-1})}^2_{H^{s+1}(R)}d\rho +K(\eta_0)\int_0^T\norm{(u^k-u^{k-1})}^2_{H^{s-1/2}(R)}d\rho\\
\leq \eta_0\int_0^T\norm{(u^k-u^{k-1})}^2_{H^{s+1}(R)}d\rho +K(\eta_0)T\sup_{t\in [0,T]}\norm{(u^k-u^{k-1})}^2_{H^{s-1/2}(R)}d\rho
\end{multline}
Combining (\ref{4.67})-(\ref{4.70}) into (\ref{4.64}) leads to
\begin{multline}\label{4.71}
\norm{u^{k+1}-u^k}^2_{H^{s-\frac{1}{2}}(R)}+|\beta|\int_0^t \norm{\partial^{s+1}_\xi (u^{k+1}-u^k)}^2_{L^2(R)} d\rho\\
\leq K_1e^{K T}C(M_0)\eta_0\int_0^t \norm{\partial^{s+1}_\xi (u^{k}-u^{k-1})}^2_{L^2(R)} d\rho\\
+K_1\left(8M_0^2e^{K T}C(M_0)\eta+ 8M_0^2e^{K T}C(M_0)K(\eta)T\right)\norm{(u^k-u^{k-1})}^2_{H^{s-1/2}(R)}\\
+K_1\left(8M_0^2e^{K T}C(M_0)T+e^{K T}C(M_0)K(\eta_0)T\right)
 \norm{(u^k-u^{k-1})}^2_{H^{s-1/2}(R)}
\end{multline}
Now choose $$\eta_0=\frac{1}{2e^{K T}C(M_0)K_1\beta}, ~~\eta=\frac{1}{32M_0^2e^{K T}C(M_0)K_1},$$
 then choose $$T\leq \frac{1}{4e^{K T}C(M_0)\left(K(\eta_0)+8M_0^2+8K(\eta)M_0^2\right)K_1},$$
 (\ref{4.71}) leads to
\begin{multline}\label{4.72}
\sup_{t\in [0,T]}\norm{u^{k+1}-u^k}^2_{H^{s-\frac{1}{2}}(R)}+|\beta|\int_0^T \norm{\partial^{s+1}_\xi (u^{k+1}-u^k)}^2_{L^2(R)} d\rho\\
\leq \frac{1}{2}\beta\int_0^T\int \norm{u^{k}-u^{k-1}}_{H^{s+1}(R)}^2  d\rho+\frac{1}{2}\sup_{t\in [0,T]} \norm{(u^k-u^{k-1})}^2_{H^{s-1/2}(R)}, 
\end{multline}

(\ref{4.72}) implies that $u^k$ is a Cauchy sequence in $C([0,T],H^{s-\frac{1}{2}}(R))$\\
$\cap L^2([0,T], H^{s+1}(R))$ and converges to a limit $u\in C([0,T],H^{s-\frac{1}{2}}(R))$\\
$\cap L^2([0,T], H^{s+1}(R))$. $ u$ is a solution of (\ref{4.1}), and equation (\ref{4.1}) implies $u\in C^{1}([0,T],H^{s-7/2}(R))$. Since $u^k$ is  bounded in 
$C([0,T],H^{s+\frac{1}{2}}(R))\cap L^2([0,T], H^{s+2}(R))$, we have $u\in L^\infty([0,T],H^{s+\frac{1}{2}}(R))\cap L^2([0,T], H^{s+2}(R)).$ Lemma 3.3 and (\ref{4.1}) imply that $u\in C([0,T],H^{s+\frac{1}{2}}(R))\cap C^{1}([0,T],H^{s-5/2}(R))$. \\

Now we prove that the solution map continuously depend on the initial data. 
Let $u_0$ and $ \tilde{u}_0$ be in $ H^{s+\frac{1}{2}}$ and $\norm{u_0}_{H^{s+\frac{1}{2}} }\le M_1$ ,$\norm{\tilde{u}_0}_{H^{s+\frac{1}{2}}}\leq M_1$ . Let $u$ and $\tilde{u}$ be solutions with initial data $u_0$ and $\tilde{u}_0$ respectively. We want to show that
\begin{equation}\label{4.73}
\sup_{t\in [0,T]}\norm{u-\tilde{u}}^2_{H^{s+\frac{1}{2}} }+\beta\int_0^T \norm{u-\tilde{u}}^2_{H^{s+2} }d\rho \leq C(M_1)\norm{u_0-\tilde{u}_0}_{H^{s+\frac{1}{2}}}.
\end{equation}

Let $w=u-\tilde{u}$, then $ w$ satisfies
\begin{equation}\label{4.74}
\begin{split}
&\partial_t w+ B[u]\tau H(\partial_\xi^3 w )
=(\mathcal{N}[u]-\mathcal{N}[\tilde{u}])+(B[\tilde u]-B[u])H[\partial_\xi^3 \tilde u] \\
&w|_{t=0}=u_0-u_0^\eta.
\end{split}
\end{equation}
 Applying Lemma 3.3 with $j=1,\epsilon=0$ to (\ref{4.74}) we obtain
\begin{multline}\label{4.75}
\norm{w}^2_{H^{s-\frac{1}{2}} }+|\beta|\int_0^t \norm{w}^2_{H^{s+1} } d\rho
\leq C\norm{u_0-\tilde{u}_0}_{H^{s-\frac{1}{2}}}\\
+K_1e^{K T} \int_0^t\norm{\mathcal{N}[u]-\mathcal{N}[\tilde u]}^2_{H^{s-2}} d\rho
+K_1e^{K T} \int_0^t\norm{(B[\tilde u]-B[u])H[\partial_\xi^3 \tilde u]}^2_{H^{s-2}} d\rho;
\end{multline}
Using the same steps as (\ref{4.65})-(\ref{4.71}), we obtain
\begin{multline}\label{4.76}
K_1e^{K T} \int_0^t\norm{\mathcal{N}[u]-\mathcal{N}[\tilde u]}^2_{H^{s-2}} d\rho
+K_1e^{K T} \int_0^t\norm{(B[\tilde u]-B[u])H[\partial_\xi^3 \tilde u]}^2_{H^{s-2}} d\rho\\
\leq \frac{1}{2}\left[\norm{w}^2_{H^{s-\frac{1}{2}} }+|\beta|\int_0^t \norm{w}^2_{H^{s+1} } d\rho\right]
\end{multline}
Combining (\ref{4.75}) and (\ref{4.76}), we obtain
\begin{equation}\label{4.77}
\norm{w}^2_{H^{s-\frac{1}{2}} }+|\beta|\int_0^t \norm{w}^2_{H^{s+1} } d\rho
\leq 2C(M_1)\norm{u_0-\tilde{u}_0}_{H^{s-\frac{1}{2}}}
\end{equation}
Now we apply Lemma 3.3 with $j=0,\epsilon=0$ to (\ref{4.74}) to obtain
\begin{multline}\label{4.78}
\norm{w}^2_{H^{s+\frac{1}{2}} }+|\beta|\int_0^t \norm{w}^2_{H^{s+2} } d\rho
\leq C\norm{u_0-\tilde{u}_0}_{H^{s+\frac{1}{2}}}\\
+K_1e^{K T} \int_0^t\norm{\mathcal{N}[u]-\mathcal{N}[\tilde u]}^2_{H^{s-1}} d\rho
+K_1e^{K T} \int_0^t\norm{(B[\tilde u]-B[u])H[\partial_\xi^3 \tilde u]}^2_{H^{s-1}} d\rho
\end{multline}
Note that
\begin{multline}\label{4.79}
\int_0^t\norm{(B[\tilde u]-B[u])H(\partial_\xi^3\tilde u)}^2_{H^{s-1}} d\rho\\
\leq \int_0^t\norm{(B[\tilde u]-B[u])}^2_{H^{s-1}}\norm{\partial_\xi^3\tilde u)}^2_{H^{s-1}} d\rho\\
\leq C(M_1)\sup_{t\in [0,T]} \norm{\tilde u-u}^2_{H^{s-1}(R)}\int_0^T\norm{\tilde u}^2_{H^{s+2}(R)}
\end{multline}
From Lemma 4.11 and (\ref{4.77}), we have
\begin{multline}\label{4.80}
\sup_{t\in [0,T]} \norm{\tilde u-u}^2_{H^{s-1}(R)}\leq 2C(M_1)\norm{u_0-\tilde{u}_0}_{H^{s-\frac{1}{2}}},~~\int_0^T\norm{\tilde u}^2_{H^{s+2}(R)}\leq 8M_1^2
\end{multline}
So (\ref{4.79}) leads to
\begin{multline}\label{4.81}
\int_0^t\norm{(B[\tilde u]-B[u])H(\partial_\xi^3\tilde u)}^2_{H^{s-1}} d\rho \\
\leq C(M_1)\norm{u_0-\tilde{u}_0}_{H^{s-\frac{1}{2}}}\leq C(M_1)\norm{u_0-\tilde{u}_0}_{H^{s+\frac{1}{2}}},
\end{multline}
Using Lemma 4.11 to obtain
\begin{multline}\label{4.82}
\norm{\mathcal{N}[u]-\mathcal{N}[\tilde u]}^2_{H^{s-1}}\\
\leq C(M_1)\norm{(u-\tilde u)}^2_{H^{s+1}(R)}+C(M_1)\norm{(u-\tilde u)}^2_{H^{3/2}(R)}\int(\partial_\xi^{s+1}\tilde u)^2d\xi
\end{multline}
Note that by Lemma 4.12 and (\ref{4.77})
\begin{multline}\label{4.83}
\sup_{t\in [0,T]}\norm{(u-\tilde u)}^2_{H^{3/2}(R)}\int_0^T\norm{\tilde u}^2_{H^{s+1}(R)}\leq ~C(M_1)\norm{u_0-\tilde{u}_0}_{H^{s+\frac{1}{2}}}~~~,
\end{multline}
and Using Lemma 2.6 we have
\begin{multline}\label{4.84}
\int_0^T\norm{(u-\tilde u)}^2_{H^{s+1}(R)}d\rho\\
\leq \eta_0\int_0^T\norm{u-\tilde u}^2_{H^{s+2}(R)}d\rho +K(\eta_0)\int_0^T\norm{u-\tilde u}^2_{H^{s+1/2}(R)}d\rho\\
\leq \eta_0 \int_0^T\norm{u-\tilde u}^2_{H^{s+2}(R)}d\rho+K(\eta_0)T\sup_{t\in [0,T]}\norm{u-\tilde u}^2_{H^{s+1/2}(R)}d\rho
\end{multline}
so (\ref{4.82}) leads to
\begin{multline}\label{4.85}
\int_0^T\norm{\mathcal{N}[u]-\mathcal{N}[\tilde u]}^2_{H^{s-1}}
\leq C(M_1)\eta_0 \int_0^T\norm{u-\tilde u}^2_{H^{s+2}(R)}d\rho\\
~~~+C(M_1)K(\eta_0)T\sup_{t\in [0,T]}\norm{u-\tilde u}^2_{H^{s+1/2}(R)}d\rho
+~C(M_1)\norm{u_0-\tilde{u}_0}_{H^{s+\frac{1}{2}}}~,
\end{multline}
choose 
$$\eta_0=\frac{1}{2e^{KT}C(M_1)\beta K_1},~~T\leq \frac{1}{2e^{KT}C(M_1)K(\eta_0)K_1};$$
and combining (\ref{4.78})-(\ref{4.85}), we obtain
\begin{multline}\label{4.86}
\norm{w}^2_{H^{s+\frac{1}{2}} }+|\beta|\int_0^t \norm{w}^2_{H^{s+2} } d\rho
\leq C(M_1)\norm{u_0-\tilde{u}_0}_{H^{s+\frac{1}{2}}}\\
+\frac{1}{2}\left[\norm{w}^2_{H^{s+\frac{1}{2}} }+|\beta|\int_0^t \norm{w}^2_{H^{s+2} } d\rho\right]
\end{multline}
which leads to (\ref{4.73}).\\
We completed the proof of Theorem 1.1.\\

{\bf Acknowledgment}: Part of this work was done when the author was in residence at MSRI in Berkeley in the Spring of 2011, the author was grateful for its support and hospitality. The author also thanks Professor Sijue Wu for very helpful discussions.


\begin{thebibliography}{99}
\bibitem{Luckhaus} S. Luckhaus, Solutions for the two-phase Stefan problem with the Gibbs-Thompson Law for the melting temperature, Euro. J. Appl. Math. {\bf 1} (1990), 101-111.
\end{thebibliography}
\end{document}